\documentclass[11pt,reqno,a4paper,twoside]{article}

\usepackage{latexsym,amsfonts,amsmath,amsthm,amstext,amscd,amssymb,euscript,wasysym}
\usepackage[hypertex]{hyperref}
\usepackage[english]{babel}
\usepackage{rotating}
\usepackage{multirow}
\usepackage{array}
\usepackage{enumitem}

\usepackage{graphicx}
\usepackage{subfigure}

\usepackage{epstopdf}

\def\k{k^{J,h}}
\def\l{l_{t,\alpha}}

\makeatletter
\@addtoreset{equation}{section}
\makeatother

\makeatletter
\@addtoreset{enunciato}{section}
\makeatother

\newcounter{enunciato}[section]

\newtheorem{ittheorem}{Theorem}
\newtheorem{itlemma}{Lemma}
\newtheorem{itproposition}{Proposition}
\newtheorem{itdefinition}{Definition}
\newtheorem{itcorollary}{Corollary}

\newenvironment{theorem}{\addtocounter{enunciato}{1}
\begin{ittheorem}}{\end{ittheorem}}

\newenvironment{lemma}{\addtocounter{enunciato}{1}
\begin{itlemma}}{\end{itlemma}}

\newenvironment{proposition}{\addtocounter{enunciato}{1}
\begin{itproposition}}{\end{itproposition}}

\newenvironment{definition}{\addtocounter{enunciato}{1}
\begin{itdefinition}}{\end{itdefinition}}
 
\newenvironment{corollary}{\addtocounter{enunciato}{1}
\begin{itcorollary}}{\end{itcorollary}}

\newcommand{\be}{\begin{equation}}
\newcommand{\ee}{\end{equation}}

\newcommand{\bt}[1]{\begin{theorem}\label{#1}}
\newcommand{\et}{\end{theorem}}

\newcommand{\bl}[1]{\begin{lemma}\label{#1}}
\newcommand{\el}{\end{lemma}}

\newcommand{\bp}[1]{\begin{proposition}\label{#1}}
\newcommand{\ep}{\end{proposition}}

\newcommand{\bd}[1]{\begin{definition}\label{#1}}
\newcommand{\ed}{\end{definition}}

\newcommand{\bpr}{\begin{proof}}
\newcommand{\epr}{\end{proof}}

\renewcommand{\phi}{\varphi}
\renewcommand{\subset}{\subseteq}
\renewcommand{\emptyset}{\varnothing}

\newcommand{\R}{\mathbb R}
\newcommand{\N}{\mathbb N}

\newcommand{\cD}{\mathcal D}

\newcommand{\cN}{\mathcal N}

\DeclareMathOperator{\csch}{csch}
\DeclareMathOperator{\acoth}{arccoth}
\DeclareMathOperator{\atanh}{arctanh}

\newcommand{\pd}[2]{\frac{\partial #1}{\partial #2}}
\newcommand{\card}[1]{\left|#1\right|}

\begin{document}

\title{Variational description of Gibbs-non-Gibbs\\ 
dynamical transitions for the Curie-Weiss model}

\author{
\renewcommand{\thefootnote}{\arabic{footnote}}
R. Fern\'andez
\footnotemark[1]
\\
\renewcommand{\thefootnote}{\arabic{footnote}}
F. den Hollander
\footnotemark[2]
\\
\renewcommand{\thefootnote}{\arabic{footnote}}
J. Mart\'inez
\footnotemark[3]
}

\footnotetext[1]{
Department of Mathematics, Utrecht University, P.O.\ Box 80010, 3508 TA Utrecht, 
The Netherlands, {\sl R.Fernandez1@uu.nl}
}
\footnotetext[2]{
Mathematical Institute, Leiden University, P.O.\ Box 9512, 2300 RA, Leiden, The 
Netherlands, \newline {\sl denholla@math.leidenuniv.nl}
}
\footnotetext[3]{
Mathematical Institute, Leiden University, P.O.\ Box 9512, 2300 RA, Leiden, The 
Netherlands, \newline {\sl martinez@math.leidenuniv.nl}
}

\maketitle

\begin{abstract}
We perform a detailed study of  
Gibbs-non-Gibbs transitions for the Curie-Weiss model subject to independent spin-flip 
dynamics (``infinite-temperature'' dynamics). We show that, in this setup, the program outlined in van Enter, Fern\'andez, den Hollander 
and Redig~\cite{vEFedHoRe10} can be fully completed, namely that Gibbs-non-Gibbs transitions are \emph{equivalent} to bifurcations in the set 
of global minima of the large-deviation rate function for the trajectories of the 
magnetization \emph{conditioned} on their endpoint. As a consequence, we show that 
the time-evolved model is non-Gibbs if and only if this set is not a singleton for 
\emph{some} value of the final magnetization. A detailed description of the possible 
scenarios of bifurcation is given, leading to a full characterization of passages from Gibbs to 
non-Gibbs ---and vice versa--- with sharp transition times (under the dynamics Gibbsianness can be lost 
and can be recovered). 

Our analysis expands the work of Ermolaev and K\"ulske~\cite{ErKu10}
who considered zero magnetic field and finite-temperature spin-flip 
dynamics. We consider both zero and non-zero magnetic field but restricted to 
infinite-temperature spin-flip dynamics. Our results reveal an interesting dependence 
on the interaction parameters, including the presence of forbidden regions for the 
optimal trajectories and the possible occurrence of overshoots and undershoots in 
the optimal trajectories. The numerical plots provided are obtained with the 
help of {\small MATHEMATICA}.

{\it Acknowledgment.} FdH is supported by ERC Advanced Grant VARIS-267356. JM is supported
by Erasmus Mundus scholarship BAPE-2009-1669. The authors are grateful to A.\ van 
Enter, V.\ Ermolaev, C.\ K\"ulske, A.\ Opoku and F.\ Redig for discussions.
\end{abstract}

\newpage


\section{Introduction and main results}
\label{S1}

Section~\ref{S1.1} provides background and motivation, Section~\ref{S1.2} a preview 
of the main results. Section~\ref{S1.3} introduces the Curie-Weiss model and the key 
questions to be explored. Section~\ref{S1.4} recalls a few facts from large-deviation 
theory for trajectories of the magnetization in the Curie-Weiss model subjected to 
infinite-temperature spin-flip dynamics and provides the link with the specification 
kernel of the time-evolved measure when it is Gibbs. Section~\ref{S1.5} states the 
main results and illustrates these results with numerical pictures. The pictures are 
made with {\small MATHEMATICA}, based on analytical expressions appearing in the text. 
Proofs are given in Sections~\ref{S2} and \ref{S3}. Section~\ref{S1} takes up half of 
the paper.


\subsection{Background and motivation}
\label{S1.1}

Dynamical Gibbs-non-Gibbs transitions represent a relatively novel and surprising 
phenome\-non. The setup is simple: an initial Gibbsian state (e.g.\ a collection of
interacting Ising spins) is subjected to a stochastic dynamics (e.g.\ a Glauber 
spin-flip dynamics) at a temperature that is \emph{different} from that of the 
initial state. For many combinations of initial and dynamical temperature, the 
time-evolved state is observed to become non-Gibbs after a finite time. Such a 
state cannot be described by any absolutely summable Hamiltonian and therefore 
\emph{lacks a well-defined notion of temperature}.   

The phenomenon was originally discovered by van Enter, Fern\'andez, den Hollander
and Redig~\cite{vEFedHoRe02} for \emph{heating dynamics}, in which a low-temperature 
Ising model is subjected to an infinite-temperature dynamics (independent spin-flips) 
or a high-temperature dynamics (weakly-dependent spin-flips). The state remains 
Gibbs for short times, but  becomes non-Gibbs after a finite time. Remarkably,
heating in this case does not lead to a succession of states with increasing 
temperature, but to states where the notion of temperature is lost altogether. 
Furthermore, it turned out that there is a difference depending on whether the 
initial Ising model has zero or non-zero magnetic field. In the former case, 
non-Gibbsianness once lost is never recovered, while in the latter case Gibbsianness 
is recovered at a later time. 

This initial work triggered a decade of developments that led to general results 
on Gibbsianness for small times (Le Ny and Redig~\cite{LNRe02}, Dereudre and 
Roelly~\cite{DeRo05}), loss and recovery of Gibbsianness for discrete spins (van 
Enter, K\"ulske, Opoku and Ruszel~\cite{KuOp08a,vEnRu08,vEnRu10,Opthesis09,vEnKuOpRu}, 
Redig, Roelly and Ruszel~\cite{ReRoRu10}), and loss and recovery of Gibbsianness 
for continuous spins (K\"ulske and Redig~\cite{KuRe06}, Van Enter and 
Ruszel~\cite{vEnRu08,vEnRu10}). A particularly fruitful research direction was 
initiated by K\"ulske and Le Ny~\cite{KuLeNy07}, who showed that Gibbs-non-Gibbs 
transitions can also be defined naturally for \emph{mean-field} models, such as 
the Curie-Weiss model. Precise results are available for the latter, including 
sharpness of the transition times and an explicit characterization of the 
conditional magnetizations leading to non-Gibbsianness (K\"ulske and 
Opoku~\cite{KuOp08b}, Ermolaev and K\"ulske~\cite{ErKu10}). In particular, the 
work in \cite{ErKu10} shows that in the mean-field setting Gibbs-non-Gibbs 
transitions occur for all initial temperatures below criticality, both for 
\emph{cooling dynamics} and for \emph{heating dynamics}.  

The ubiquitousness of the Gibbs-non-Gibbs phenomenon calls for a better 
understanding of its causes and consequences. Unfortunately, the mathematical 
approach used in most references is opaque on the intuitive level. Generically, 
non-Gibsianness is proved by looking at the evolving system at two times, 
the inital and the final time, and applying techniques from equilibrium 
statistical mechanics. This is an indirect approach that does not illuminate 
the relation between the Gibbs-non-Gibbs phenomenon and the dynamical effects 
responsible for its occurrence.  This unsatisfactory situation was addressed 
in Enter, Fern\'andez, den Hollander and Redig~\cite{vEFedHoRe10}, where possible 
dynamical mechanisms were proposed and a \emph{program} was put forward to 
develop a theory of Gibbs-non-Gibbs transitions on \emph{purely dynamical 
grounds}. The present paper shows that this program can be fully carried out 
for the Curie-Weiss model subject to an infinite-temperature dynamics. 

In the mean-field scenario, the key object is the time-evolved single-spin average 
\emph{conditional} on the final empirical magnetization. Non-Gibbsianness 
corresponds to a discontinuous dependence of this average on the final 
magnetization. The discontinuity points are called \emph{bad magnetizations} 
(see Definition~\ref{goodpoint} below). Dynamically, such discontinuities are 
expected to arise whenever there is more than one possible trajectory \emph{compatible} 
with the bad magnetization at the end. Indeed, this expectation is confirmed and 
exploited in the sequel. The actual conditional trajectories are those minimizing 
the large-deviation rate function on the space of trajectories of magnetizations. 
The time-evolved measure remains Gibbsian whenever there is a single minimizing 
trajectory for every final magnetization, in which case the specification kernel 
can be computed explicitly (see Proposition~\ref{Thetheorem} below). In contrast, 
if there are multiple optimal trajectories, then the choice of trajectory can be 
decided by an infinitesimal perturbation of the final magnetization, and this is 
responsible for non-Gibbsianness.


\subsection{Preview of the main results}
\label{S1.2}

In the present paper we study in detail the large-deviation rate function for the 
trajectory of the magnetization in the Curie-Weiss model with pair potential $J>0$ 
and magnetic field $h \in \R$ (see \eqref{Hdef} below). We exploit the fact that, 
due to the mean-field character of the interaction, this rate function can be 
expressed as a function of the initial and the final magnetization only (see
Proposition~\ref{LDP} below), i.e., the trajectories are uniquely determined 
by the magnetizations at the beginning and at the end (see Corollary~\ref{CIrelation} 
and Proposition~\ref{Thereduction} below). Here is a summary of the main results
(see Fig.~\ref{fig-crossovertimes}):
\begin{enumerate}
\item 
If $0 < J \leq 1$ (supercritical temperature), then the evolved state is Gibbs at all 
times. On the other hand, if $J > 1$ (subcritical temperature) there exists 
some time $\Psi_U$ at which multiple trajectories appear. The associated non-Gibbsianness persists for all later times when $h=0$ (zero magnetic field). All these features were 
already shown by Ermolaev and K\"ulske~\cite{ErKu10}.
\item 
For $h \neq 0$ there is a time $\Psi_*>\Psi_U$ at which Gibbsianness is restored for 
all later times.  
\item 
There is a change in behavior at $J=\tfrac32$. For $1 < J\leq \tfrac32$:
\begin{enumerate}
\item 
If $h=0$, then only the zero magnetization is bad for $t > \Psi_c$.
\item 
If $h>0$ ($h<0$), then there is only one bad magnetization for $\Psi_U < t \leq \Psi_*$. 
This bad magnetization changes with $t$ but is always strictly negative (strictly 
positive).
\end{enumerate}
For $J>\tfrac32$:
\begin{enumerate}
\item 
If $h=0$, then there is a time $\Psi_c>\Psi_U$ such that for $\Psi_U<t<\Psi_c$ there 
are two non-zero bad magnetizations (equal in absolute value but with opposite signs), 
while for $t \geq \Psi_c$ only the zero magnetization is bad.
\item 
If $h \neq 0$ and small enough, then there are two times $\Psi_T>\Psi_L$ between
$\Psi_U$ and $\Psi_*$ such that for $\Psi_U < t \leq \Psi_L$ and $\Psi_T \leq t 
\leq \Psi_*$ only one bad magnetization occurs, while for $\Psi_L < t < \Psi_T$ 
two bad magnetizations occur.
\end{enumerate}
\end{enumerate}
  
\begin{figure}[htbp]
\setlength{\unitlength}{0.4cm}
\begin{center}
\begin{picture}(10,10)(3,0)
\put(0,0){\line(18,0){18}}
\put(0,3){\line(18,0){18}}
\put(-3,2.7){$h=0$}
\put(-3,-.3){$h \neq 0$}
\put(2,0){\circle*{.35}}
\put(2,3){\circle*{.35}}
\put(4,0){\circle*{.35}}
\put(10,0){\circle*{.35}}
\put(15,0){\circle*{.35}}
\put(9,3){\circle*{.35}}
\put(1.6,-1.3){$\Psi_U$}
\put(14.6,-1.3){$\Psi_*$}
\put(3.6,-1.3){$\Psi_L$}
\put(9.6,-1.3){$\Psi_T$}
\put(1.6,3.8){$\Psi_U$}
\put(8.6,3.8){$\Psi_c$}
\end{picture}
\end{center}
\caption{Crossover times for $h=0$ and $h \neq 0$ when $J>\tfrac32$.}
\label{fig-crossovertimes}
\end{figure}
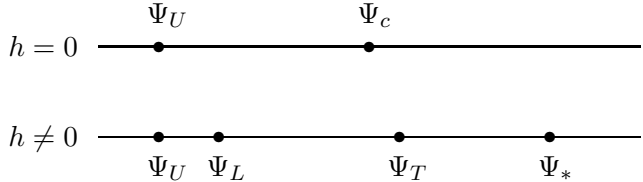

All the crossover times depend on $J,h$ and are strictly positive and finite.
Our analysis gives a detailed picture of the optimal trajectories for different 
$J,h$ and different conditional magnetizations. Among the novel features we mention:
\begin{enumerate}
\item[(1)] 
Presence of \emph{forbidden regions} that cannot be crossed by any optimal 
trajectory. The boundary of these regions is given by the multiple optimal 
trajectories when bifurcation sets in.   The forbidden regions were predicted in~\cite{vEFedHoRe10} and first found, for $h=0$, by Ermolaev and K\"ulske~\cite{ErKu10}. 
\item[(2)] 
Existence of \emph{overshoots} and \emph{undershoots} for optimal trajectories
for $h \neq 0$.
\item[(3)] 
Classification of the bad magnetizations leading to multiple optimal trajectories. 
These bad magnetizations depend on $J,h$ and change with time.
\end{enumerate}


\subsection{The model}
\label{S1.3}


\subsubsection{Hamiltonian and dynamics}
\label{S1.3.1}

The Curie-Weiss model consists of $N$ Ising spins, labelled $i=1,\ldots,N$ with
$N \in \N$. The spins interact through a mean-field Hamiltonian ---that is, a Hamiltonian involving no geometry and no sense of neighborhood, in which each spin interacts equally with all other spins---.   
The Curie-Weiss Hamiltonian is
\begin{equation}
\label{Hdef}
H^N(\sigma) := -\tfrac{J}{2N} \sum\limits_{i,j=1}^{N} \sigma_i \sigma_j 
- h \sum\limits_{i=1}^{N} \sigma_i, \qquad \sigma\in\Omega_N,
\end{equation}
where $J>0$ is the (ferromagnetic) pair potential, $h\in\R$ is the (external) 
magnetic field, $\Omega_N:=\{-1,+1\}^N$ is the spin configuration space, and 
$\sigma:=(\sigma_i)_{i=1}^N$ is the spin configuration. The Gibbs measure associated 
with $H^N$ is 
\begin{equation}
\label{CWGibbs}
\mu^N(\sigma) := \frac{e^{-H^N(\sigma)}}{Z^N}, \qquad \sigma\in\Omega_N,
\end{equation}
with $Z^N$ the normalizing partition sum. 

We allow this model to evolve according to an  
\emph{independent spin-flip} dynamics, that is, a dynamics defined by the  generator $L_N$ given by (see 
Liggett~\cite{Li85} for more background)
\begin{equation}
(L_N f)(\sigma)
:= \sum\limits_{i=1}^{N} [f(\sigma^i)-f(\sigma)], \quad  f\colon\,\Omega_N\to\R,
\end{equation}
where $\sigma^i$ denotes the configuration obtained from $\sigma$ by flipping the 
spin with label $i$.  The resulting random variables $\sigma(t):=(\sigma_i(t))_{i=1}^N$ 
constitute a continuous-time Markov chain on $\Omega_N$.   
We write $\mu^N_t$ to denote the measure on $\Omega_N$ at 
time $t$ when the initial measure is $\mu^N$ and abbreviate $\mu_t:=
(\mu^N_t)_{N\in\N}$.


\subsubsection{Empirical magnetization}
\label{S1.3.2}

To emphasize its mean-field character, it is convenient to write the Hamiltonian (\ref{Hdef}) in the form
\begin{equation}
\label{Hdefmf}
H^N(\sigma) = N\bar H(m_N(\sigma))
\end{equation}
 where
 \be
\bar{H}(x) := -\tfrac12J x^2-h x, \qquad x\in\R.
\ee  
and 
 \begin{equation}
\label{empmag}
m_N(\sigma) := \frac{1}{N} \sum_{i=1}^{N}\sigma_i
\end{equation}
is the empirical magnetization of $\sigma\in\Omega_N$, which takes values in the set $\mathcal{M}_N := \{-1,-1+2N^{-1},\ldots,+1-2N^{-1},+1\}$. 
The Gibbs measure on $\Omega_N$ induces a Gibbs measure on $\mathcal{M}_N$
given by
\begin{equation}
\label{CW}
\bar\mu^N(m) := \binom{N}{\tfrac{1+m}{2}\,N}\,
\frac{e^{-N\bar{H}(m)}}{\bar{Z}^N}, \qquad m \in \mathcal{M}_N,
\end{equation}
where $\bar{Z}^N$ is the normalizing partition sum.

The independent (infinite-temperature) dynamics has the simplifying feature of preserving the mean-field character of the model.  In fact, the dynamics on $\Omega_N$ induces a 
dynamics on $\mathcal{M}_N$, which is a continuous-time Markov chain $(m^N_t)_{t
\geq 0}$ with generator $\bar{L}_N$ given by
\be
(\bar{L}_Nf)(m) := \frac{1+m}{2}\,N[f(m-2N^{-1})-f(m)]
+\frac{1-m}{2}\,N[f(m+2N^{-1})-f(m)]\;, 
\ee
for \quad $f\colon\,\mathcal{M}_N \to \R$.
Adapting our previous notation we denote
$\bar\mu^N_t$ the measure on $\mathcal{M}_N$ at time $t$, and 
abbreviate $\bar\mu_t:=(\bar\mu^N_t)_{N\in\N}$. Due to permutation 
invariance, $\mu^N_t$ characterizes $\bar\mu^N_t$ and vice versa, for each $N$ and $t$. We write $P^N$ 
to denote the law of $(m^N_t)_{t\geq 0}$, which lives on the space of c\`adl\`ag 
trajectories $D_{[0,\infty)}([-1,+1])$ endowed with the Skorohod topology.


\subsubsection{Bad magnetizations}
\label{S1.3.3}

Non-Gibbsianness shows up through discontinuities with respect to boundary conditions
of finite-volume conditional probabilities.  For the Curie-Weiss model it is enough to
consider the single-spin conditional probabilities
\be
\label{specification}
\gamma^N_t(\sigma_1 \mid \alpha_{N-1}) := \mu^N_t(\sigma_1 \mid \sigma_{N-1})\;,
\ee
defined for $\sigma_1 \in \{-1,+1\}$ and $\alpha_{N-1} \in \mathcal{M}_{N-1}$, 
and any spin configuration $\sigma_{N-1}\in\Omega_{N-1}$ such that 
$m_{N-1}(\sigma_{N-1})=\alpha_{N-1}$. By permutation invariance, \eqref{specification} 
does not depend on the choice of $\sigma_{N-1}$.

The central definition for our purposes is the following.
\begin{definition}
\label{goodpoint}
{\rm (K\"ulske and Le Ny~\cite{KuLeNy07})} Fix $t \geq 0$.\\ 
(a) A magnetization $\alpha \in [-1,+1]$ is said to be good for $\mu_t$ if there exists a 
neighborhood $\mathcal{N}_\alpha$ of $\alpha$ such that
\be
\label{specificationlimit}
\gamma_t(\cdot \mid \bar{\alpha}) 
:= \lim\limits_{N\to\infty} \gamma^N_t(\cdot \mid \alpha_{N-1}),
\ee
exists for all $\bar{\alpha}$ in $\mathcal{N}_\alpha$ and all $(\alpha_N)_{N\in\N}$ such 
that $\alpha_N\in\mathcal{M}_N$ for all $N \in \N$ and $\lim_{N\to\infty} \alpha_N
=\bar{\alpha}$, and is independent of the choice of $(\alpha_N)_{N\in\N}$. The limit 
is called the specification kernel. In particular, $\bar{\alpha} \mapsto \gamma_t(\cdot 
\mid \bar{\alpha})$ is continuous at $\bar{\alpha}=\alpha$.\\
(b) A magnetization $\alpha \in [-1,+1]$ is called bad if it is not good.\\
(c) $\mu_t$ is called Gibbs if it has no bad magnetizations.
\end{definition}


\subsection{Path large deviations and link to specification kernel}
\label{S1.4}

The main point of our work is our relation between path large deviations and non-Gibbsianness.  
For the convenience of the reader, let us recall some basic large deviation results for the Curie-Weiss model. 
For background on large deviation theory, see e.g.\ den Hollander~\cite{dHo00}.


\subsubsection{Path large deviation principle}
\label{S1.4.1}

Let us recall that a family of measures $\nu^N$ on a Borel measure space satisfies a \emph{large deviation principle} with rate function $I$ and speed $N$ if the following two conditions are satisfied: 
\begin{eqnarray}
\liminf_{N\to\infty} \frac1N \log \nu^N(A) &\ge & -\inf_{x\in A} I(x) \qquad \mbox{for } A \mbox{ open}\\
\limsup_{N\to\infty} \frac1N \log \nu^N(A) &\le & -\sup_{x\in A} I(x) \qquad \mbox{for } A \mbox{ closed}
\end{eqnarray}

The proof of the following proposition is elementary and can be found in many
references. The indices $S$ and $D$ stand for \emph{static} and \emph{dynamic}.

\begin{proposition}
\label{LDP}
{\rm (Ermolaev and K\"ulske~\cite{ErKu10}, Enter, Fern\'andez, den Hollander 
and Redig~\cite{vEFedHoRe10})}\\
(i) $(\bar\mu^N)_{N\in\N}$ satisfies the large deviation principle on $[-1,+1]$ with rate 
$N$ and rate function $I_S-\inf(I_S)$ given by
\be
I_S(m) := \bar{H}(m)+\bar{I}(m), \qquad 
\bar{I}(m):=\frac{1+m}{2}\,\log(1+m)+\frac{1-m}{2}\,\log(1-m).
\ee
(ii) For every $T>0$, the restriction of $(P^N)_{N\in\N}$ to the time interval $[0,T]$
satisfies the large deviation principle on $D_{[0,T]}([-1,+1])$ with rate $N$ and rate 
function $I^T-\inf(I^T)$ given by
\be
\label{originalcost}
I^T(\phi) := I_S(\phi(0)) + I_D^T(\phi),
\ee
where
\be
\label{actionid}
I_D^T(\phi) := \left\{\begin{array}{ll}
\int_0^T L(\phi(s),\dot{\phi}(s))\,ds &\mbox{ if } \dot\phi \mbox{ exists},\\
\infty &\mbox{ otherwise},
\end{array}
\right.
\ee
is the action integral with Lagrangian
\be
\label{Lform}
L(m,\dot{m}) := -\frac12 \sqrt{4 \left(1-m^2\right)+\dot{m}^2}
+\frac12 \dot{m} \log 
\left(\frac{\sqrt{4 \left(1-m^2\right)+\dot{m}^2}+\dot{m}}{2(1-m)}\right)+1.
\ee
\end{proposition}

Let 
\be
Q^N_{t,\alpha}(m) := P^N \left(m_N(0)=m \mid m_N(t)=\alpha \right), 
\qquad m\in\mathcal{M}_N
\ee
be the conditional distribution of the 
magnetization at time 0 given that the magnetization at time $t$ is $\alpha$.
The contraction principle applied to Proposition~\ref{LDP}(ii) implies the 
following large deviation principle.

\begin{corollary}
\label{CIrelation}
For every $t \geq 0$ and $\alpha \in [-1,+1]$, $(Q^N_{t,\alpha})_{N \in \N}$ 
satisfies the large deviation principle on $[-1,+1]$ with rate $N$ and rate 
function $C_{t,\alpha}-\inf(C_{t,\alpha})$ given by
\be
\label{variational}
C_{t,\alpha}(m)
:= \inf_{\substack{\phi\colon\,\phi(0)=m, \\ \phi(t)=\alpha}} I^t(\phi).
\ee
\end{corollary}

\noindent
Note that
\be
\label{cvp}
\inf_{m \in [-1,+1]} C_{t,\alpha}(m) = \inf_{m \in [-1,+1]}\,\, 
\inf_{\substack{\phi\colon\,\phi(0)=m, \\ \phi(t)=\alpha}} 
I^t(\phi)= \inf_{\phi\colon\,\phi(t)=\alpha} I^t(\phi).
\ee


\subsubsection{Link to specification kernel}
\label{S1.4.2}

The following proposition provides the fundamental link between the specification 
kernel in \eqref{specificationlimit} and the minimizer of \eqref{cvp} when it is 
unique, and is a straightforward generalization to arbitrary magnetic field of a 
result for zero magnetic field stated and proved in Ermolaev and K\"ulske~\cite{ErKu10}.

\begin{proposition}
\label{Thetheorem}
Fix $t \geq 0$ and $\alpha \in [-1,+1]$. Suppose that \eqref{cvp} has a unique 
minimizing path $(\hat{\phi}_{t,\alpha}(s))_{0 \leq s \leq t}$. Then the specification 
kernel equals
\be
\label{minvsspecif}
\gamma_t(z \mid \alpha) =
\frac{\sum_{x\in\{-1,+1\}} e^{x[J\hat{\phi}_{t,\alpha}(0)+h]} p_t(x,z)}
{\sum_{x,y\in\{-1,+1\}} e^{x[J\hat{\phi}_{t,\alpha}(0)+h]} p_t(x,y)}, \qquad 
z \in \{-1,+1\},
\ee
where $p_t(\cdot,\cdot)$ is the transition kernel of the continuous-time Markov chain 
on $\{-1,+1\}$ jumping at rate $1$, given by $p_t(1,1)=p_t(-1,-1)=e^{-t} \cosh(t)$ and
$p_t(-1,+1)=p_t(1,-1)=e^{-t} \sinh(t)$.
\end{proposition}

\noindent
{\bf Remark:}
Note that the expression in the right-hand side of \eqref{minvsspecif} depends on 
the optimal trajectory only via its initial value $\hat{\phi}_{t,\alpha}(0)$. Thus, 
\eqref{minvsspecif} has the form
\begin{equation}
\label{spec-initialm}
\gamma_t(z \mid \alpha)=\Gamma_t(z,J \hat{\phi}_{t,\alpha}(0)+h),
\end{equation}
where $\hat{\phi}_{t,\alpha}(0)$ is the unique global minimizer of $m \mapsto 
C_{t,\alpha}(m)$ and $m \mapsto \Gamma_t(z,m)$ is continuous and strictly increasing 
(strictly decreasing) for $z=1$ ($z=-1$).


\subsubsection{Reduction}
\label{S1.4.3}

The next proposition allows us to reduce \eqref{cvp} to a one-dimensional variational 
problem. Consider the equation
\be
\label{Pequation}
k^{J,h}(m)=l_{t,\alpha}(m)
\ee
with
\be
\label{kbh}
\begin{aligned}
k^{J,h}(m) &:= a^J(m) \cosh(2h) + b^J(m) \sinh(2h),\\
l_{t,\alpha}(m) &:=m \coth(2t) - \alpha \csch(2t),
\end{aligned}
\ee
where
\be
\label{gb}
\begin{aligned} 
a^J(m) &:=\sinh(2Jm) - m \cosh(2Jm),\\
b^J(m) &:=\cosh(2Jm) - m \sinh(2Jm).
\end{aligned}
\ee

\begin{proposition}
\label{Thereduction}
Let $C_{t,\alpha}$ be as in \eqref{variational}. Then, for every $t \geq 0$ and 
$\alpha \in [-1,+1]$,
\be
\label{cost}
\begin{split}
&C_{t,\alpha}(m) = I_S(m)\\ 
&\qquad + \frac14 \left\{4t +\log\left(\frac{1-\alpha^2}{1-m^2}\right)
+ \log\left(\left[\frac{1-R-2C_1\alpha e^{-2t}}{1+R-2C_1\alpha e^{-2t}}\right] 
\left[\frac{1+R-2C_1 m}{1-R-2C_1 m}\right]\right)\right.\\
&\qquad \left. + 2 \left[\alpha\log\left(\frac{R-C_1e^{-2t}+C_2e^{2t}}{1-\alpha}\right)
- m \log\left(\frac{R-C_1+C_2}{1-m}\right)\right]\right\}
\end{split}
\ee			
with
\be \label{cost.1}
\begin{array}{lllll}
C_1 &=& C_1(t,\alpha,m) &:=& \frac{m e^{2t}-\alpha}{e^{2t}-e^{-2t}},\\[0.2cm]
C_2 &=& C_2(t,\alpha,m) &:=& \frac{\alpha-me^{-2t}}{e^{2t}-e^{-2t}},\\[0.2cm]
R &=& R(C_1,C_2) &:=& \sqrt{1-4C_1C_2}. 
\end{array}
\ee
Furthermore, the critical points of $C_{t, \alpha}$ are the solutions of \eqref{Pequation}. Hence,
\be
\label{simplification}
\inf_{\substack{\phi\colon\,\phi(t)=\alpha}} I^t(\phi)
= \min_{m \text{ solves } \eqref{Pequation}} C_{t,\alpha}(m)\;,
\ee
and the constrained minimizing trajectories are of the form
\begin{eqnarray}
\label{trajectory}
\hat{\phi}^{\hat{m}}_{t,\alpha}(s) &:=& \csch(2t) \Big\{m \sinh(2(t-s)) + \alpha \sinh(2s) \Big\}\qquad
0\le s\le t\\
\label{trajectory.00}
\hat{m} =\hat{m}(t,\alpha)&=& {\rm argmin} \Bigl[C_{t,\alpha}\bigr |_{\text{solutions of } \eqref{Pequation}}\Bigr]\;.
\end{eqnarray}
\end{proposition}

\noindent 
The identities
\be\label{eq:rr.ident}
k^{J,h}(m)\;=\; 2 \cosh^2(Jm+h) \bigl[\tanh(Jm + h)-m\bigr] + m
\ee
and
\be
\lim_{t\to\infty} l_{t,\alpha}(m) = m
\ee
imply that
in the limit $t\to\infty$ \eqref{Pequation} reduces to $\tanh(Jm + h)=m$.
This is the equation for the spontaneous magnetization of the Curie-Weiss model 
with parameters $J,h$. This equation has always at least one solution and the value
\be
m^\infty = m^\infty(J,h) 
:= \mbox{ the largest solution of the equation } \tanh(Jm + h)=m
\ee
is well known to be strictly positive if $h>0$ or if $J>1$.  In these regimes, the standard Curie-Weiss graphical argument shows that, for $m>0$,
\be
\label{eq:rr.infty}
k^{J,h}(m) \;\substack{<\\=\\>}\; m \quad \Longleftrightarrow\quad 
m \;\substack{>\\=\\<}\; m^\infty\;. 
\ee
We also remark that when $t\to 0$ the function $l_{t,\alpha}$ converges to the line defined by the equation $m=\alpha$. This implies that for short times there is a unique solution of \eqref{Pequation} and it is close to $\alpha$.

\subsection{Main results}
\label{S1.5}

In Section~\ref{S1.5.1} we state the equivalence of non-Gibbs and 
bifurcation that lies at the heart of the program outlined in \cite{vEFedHoRe10} 
(Theorem~\ref{nongibbs-bifurcation}). In Section~\ref{S1.5.2} we introduce
some notation. In Section~\ref{S1.5.3} we identify the optimal trajectories for 
$\alpha=0$, $h=0$ (Theorems~\ref{hmfzero}--\ref{Cone-Monotonicity}). In 
Section~\ref{S1.5.4} we extend this identification to $\alpha \in [-1,+1]$, $h\in\R$
(Theorem~\ref{generaltheorem}). In Section~\ref{S1.5.5} we summarize the consequences
for Gibbs versus non-Gibbs (Corollary~\ref{GNG-coro}).


\subsubsection{Equivalence of non-Gibbs and bifurcation}
\label{S1.5.1}

The following theorem proves the long suspected equivalence between 
dynamical non-Gibbsianness, i.e., discontinuity of $\alpha \mapsto \gamma_t(\cdot \mid \alpha)$ at $\alpha_0$, 
and non-uniqueness of the global minimizer of $m \mapsto C_{t,\alpha_0}(m)$, i.e.,
the occurrence of more than one possible history for the same $\alpha$.

\begin{theorem}
\label{nongibbs-bifurcation} 
$\alpha \mapsto \gamma_t(\sigma \mid \alpha)$ is continuous at $\alpha_0$ if and 
only if $\inf_{\phi\colon\,\phi(t)=\alpha_0} I^t(\phi)$ has a unique minimizing 
path or, equivalently, $\inf_{m\in [-1,+1]} C_{t,\alpha_0}(m)$ has a unique minimizing 
magnetization.
\end{theorem}


\subsubsection{Notation}
\label{S1.5.2} 

Due to relation \eqref{simplification}, our analysis focusses on the different solutions of \eqref{Pequation} obtained as $t,\alpha$ are varied.  In particular, we must determine which of them are minima of the variational problem in \eqref{cvp}. We 
write 
\be
\Delta_{t,\alpha} := \mbox{ the set of global minimizers of } C_{t,\alpha}.
\ee 
For brevity, when $\alpha$ is kept fixed and $\Delta_{t,\alpha}$ is a singleton $\{\hat{m}(t,\alpha)\}$ for each $t$, we write $\hat{m}(t)$ instead of $\hat{m}(t,\alpha)$.
When $h,\alpha=0$, by symmetry we have $\Delta_{t,0}=\{0\}$ or 
$\Delta_{t,\alpha}=\{\pm\hat{m}(t)\}$, where in the last case we denote by 
$\hat{m}(t)$ the unique positive global minimizer.
If both the initial and final magnetizations are fixed, then there is a unique minimizer that we denote as in \eqref{trajectory}.  That is,
\be
\hat{\phi}^m_{t,\alpha} := \mathop{\mbox{argmin}}\limits_{\substack{\phi\colon\,\phi(0)=m, \\ \phi(t)=\alpha}} I^t_D(\phi)
\ee 
for $m,\alpha\in[-1,+1]$.
We emphasize that, by definition, $C_{t,\alpha}(m)=I^t(\hat{\phi}^{m}_{t,\alpha})$
and $\hat{\phi}_{t,\alpha}(s)=\hat{\phi}^{\hat{m}(t,\alpha)}_{t,\alpha}(s), 
s \in [0,t]$.  In particular $\hat{m}(t, \alpha)=\hat{\phi}_{t,\alpha}(0)$.


\subsubsection{Optimal trajectories for $\alpha=0$, $h=0$}
\label{S1.5.3}

The following theorem refers to a critical time  
\be
\label{CRITICALT}
\Psi_c = \Psi_c(J):=
\begin{cases}
\tfrac12 \acoth(2J -1) &\text{if } 1 < J \leq \tfrac32,\\
t_* &\text{if } J > \tfrac32,
\end{cases}
\ee
where $t_*=t_*(J)$ is implicitly calculable: $t_*=t(m_*)$ where the function $t(m)$ is defined 
in \eqref{tm} below and $m_*=m_*(J)$ is the solution of \eqref{mjump}.

\begin{figure}[htbp]
\label{theorem1}
\begin{center}
\begin{tabular}{ccc}
\includegraphics[width=0.25\textwidth]{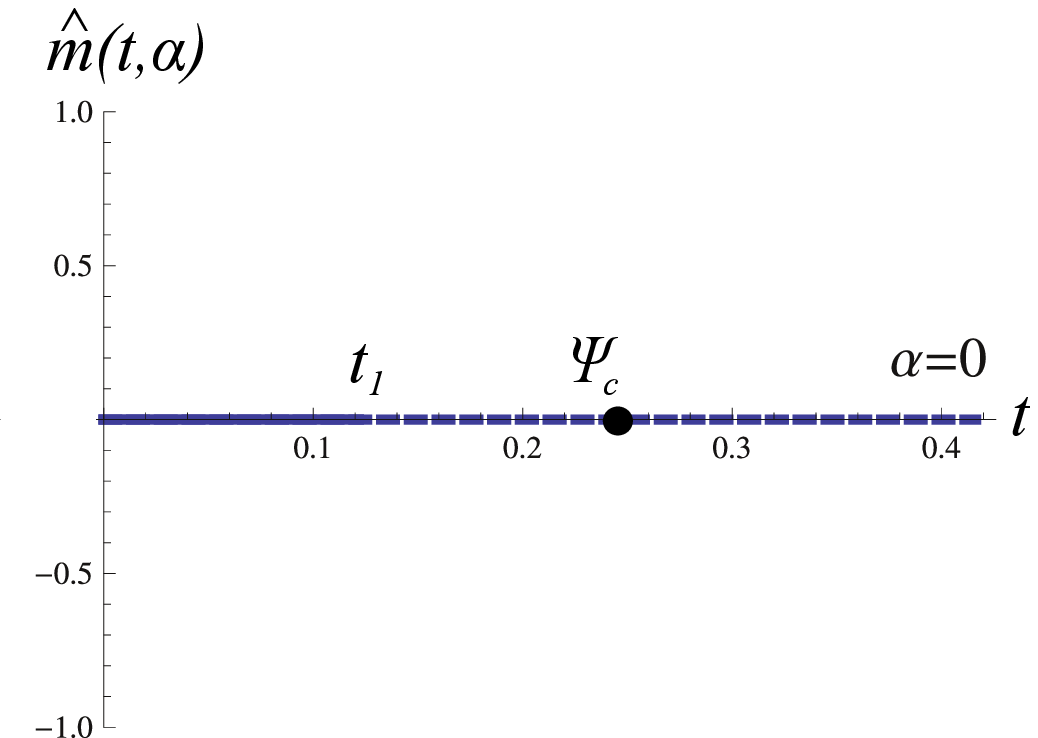} &
\includegraphics[width=0.25\textwidth]{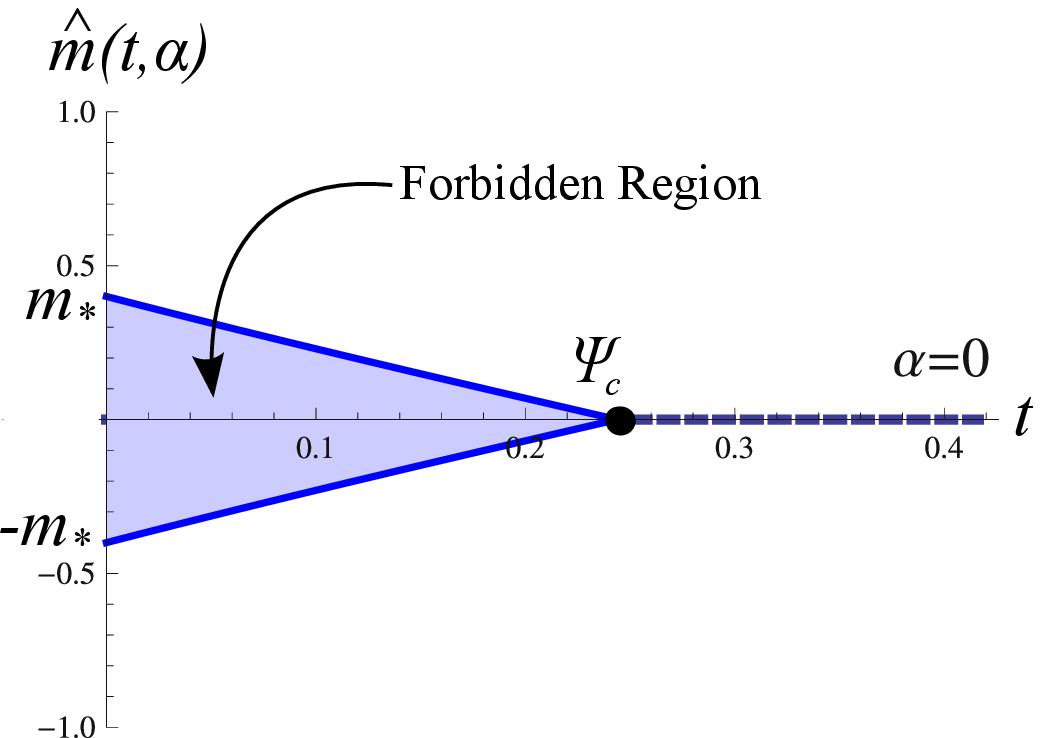} &
\includegraphics[width=0.25\textwidth]{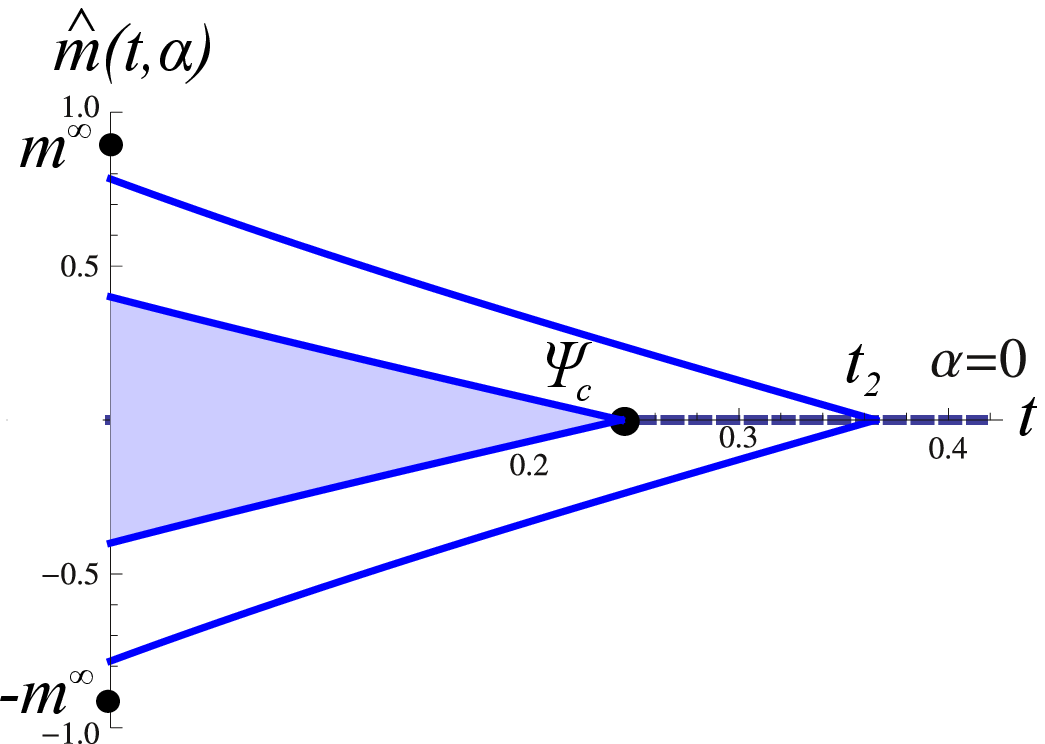} \\
\includegraphics[width=0.25\textwidth]{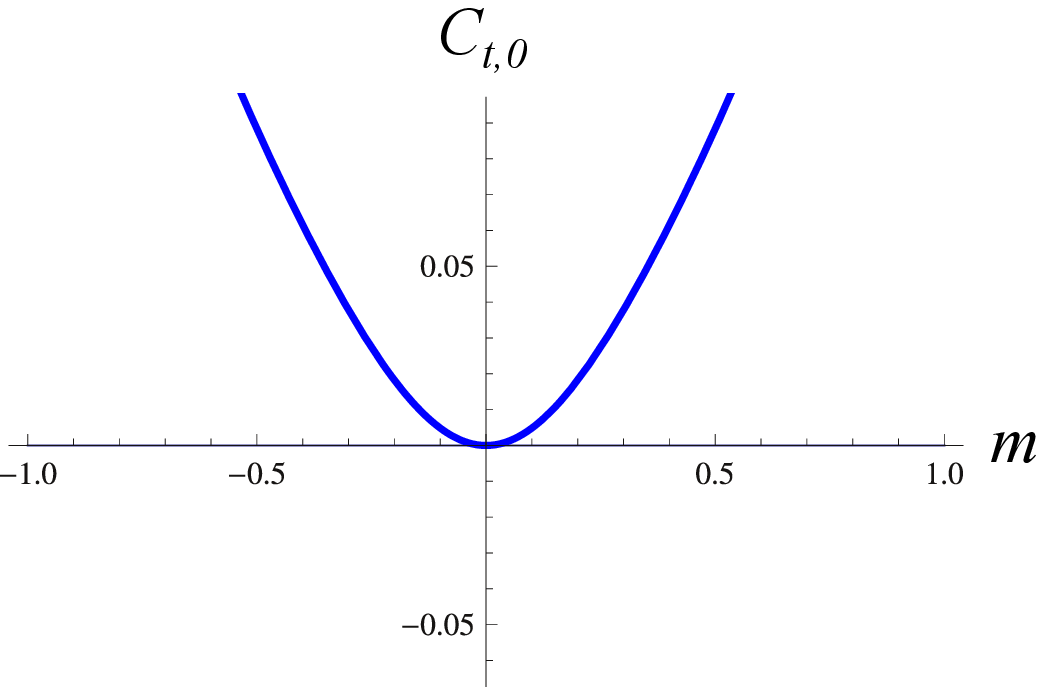} &
\includegraphics[width=0.25\textwidth]{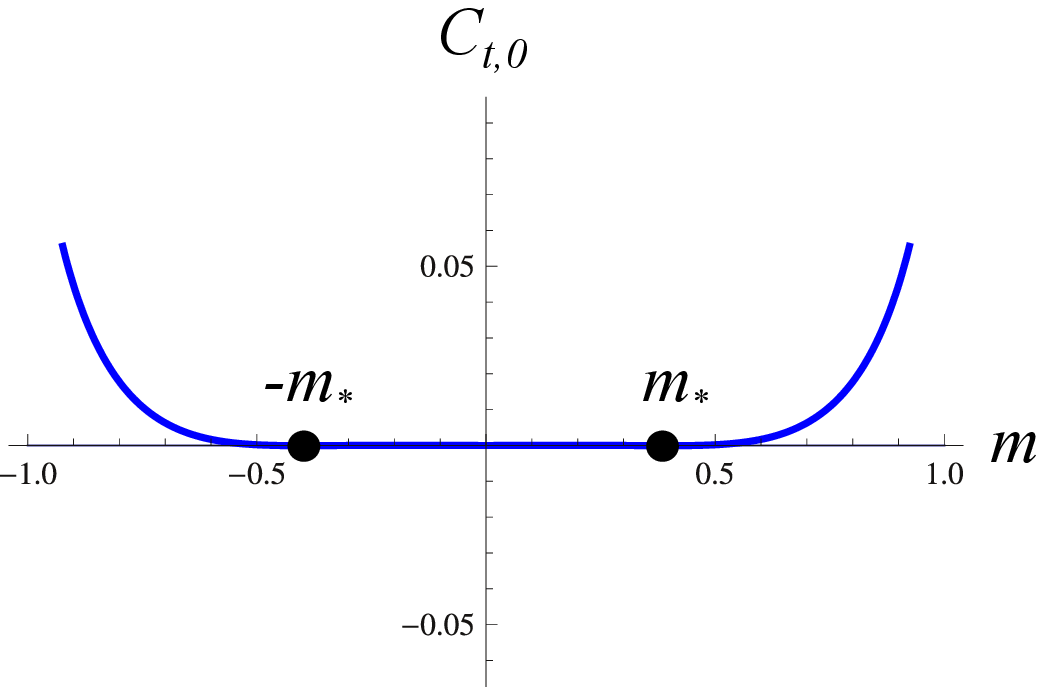} &
\includegraphics[width=0.25\textwidth]{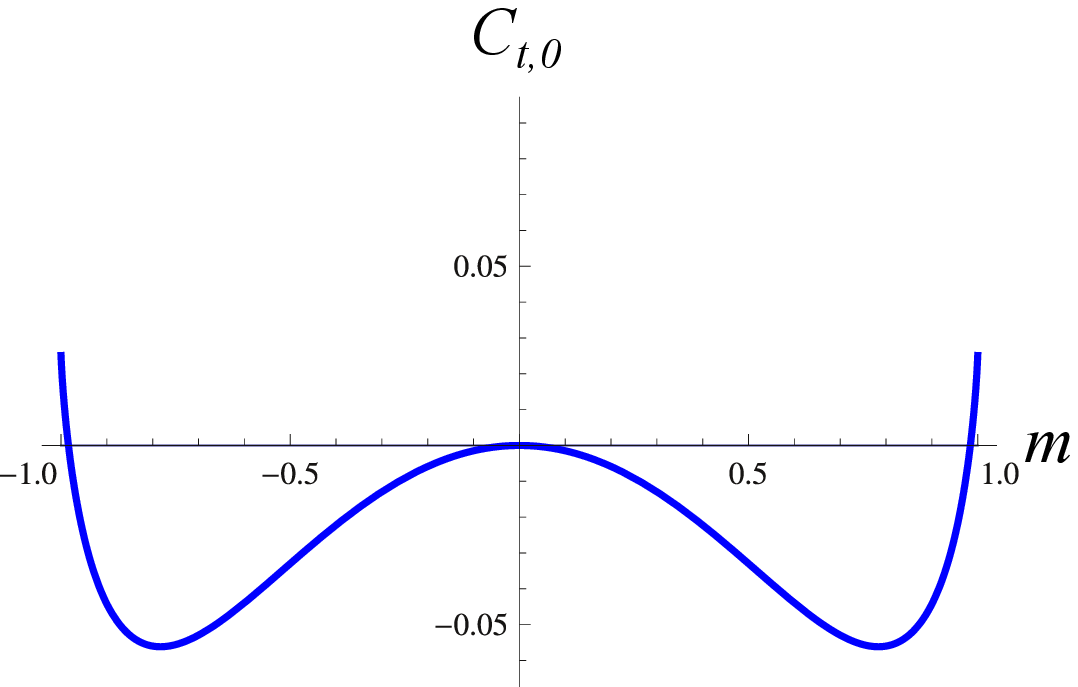} \\
&&\\
$t(=t_1)<\Psi_c$ & $t=\Psi_c$ & $t(=t_2)>\Psi_c$
\end{tabular}
\end{center}
\caption{Illustration of Theorem \ref{hmfzero}. \emph{First row:} Time evolution
of the minimizing trajectories $\pm(\hat{\phi}_{t,0}(s))_{0\leq s \leq t}$ for $t<\Psi_c$, 
$t=\Psi_c$ and $t>\Psi_c$ for an initial Curie-Weiss model with $(J,h)=(1.6,0)$ [regime (iii) in the Theorem]. The shaded cone is the forbidden region. \emph{Second row:} 
Plot of $m \mapsto C_{t,0}(m)$ for the same times and parameter values.}
\label{fig-1}
\end{figure}

\begin{theorem}
\label{hmfzero}
{\rm (See Fig.~\ref{fig-1}.)}
Consider $\alpha=0$ and $h=0$.
\begin{itemize}
\item[(i)] If $0<J\leq 1$, then 
\be
\Delta_{t,0}=\{0\}, \quad \forall\,t \geq 0.
\ee
\item[(ii)] If $1<J\leq\tfrac32$, then
\be
\Delta_{t,0} = \left\{\begin{array}{ll}
\{0\} &\text{if } 0 \leq t \leq \Psi_c,\\
\{ \pm \hat{m}(t)\} &\text{if } t > \Psi_c,
\end{array}
\right.
\ee
where $t\mapsto\hat{m}(t)$ is continuous and strictly increasing on $[\Psi_c,\infty)$ with 
$\hat{m}(\Psi_c)=0$.
\item[(iii)] If $J>\tfrac32$, then
\be
\Delta_{t,0} = \left\{\begin{array}{ll}
\{0\} &\text{if } 0 \leq t < \Psi_c,\\
\{\pm \hat{m}(t)\} &\text{if } t \geq \Psi_c,
\end{array}
\right.
\ee
where $t\mapsto\hat{m}(t)$ is continuous and strictly increasing on $[\Psi_c,\infty)$ with 
$\hat{m}(\Psi_c)=:m_*>0$.
\end{itemize}
\end{theorem}

Let $\Lambda_{t,0}(J)$ denote the cone between the trajectories $\pm\hat{\phi}_{t,0}$. 
As a consequence of the previous theorem, \emph{no} minimal trajectory conditioned in $t'$ with $t'\geq t$ can intersect the interior of this region.
Such a cone corresponds, therefore, to a \emph{forbidden region}. 
Forbidden regions grow, in a nested fashion, as the conditioning time $t$ grows.
There is, however, a distinctive difference between regimes (ii) and (iii) in the previous theorem:
In Regime (ii) the forbidden region opens up \emph{continuously} after $\Psi_c$, 
while in Regime (iii) it opens up \emph{discontinuously}. These facts are summarized in the following theorem.

\begin{theorem}
\label{Cone-Monotonicity}
Suppose that $\alpha=0$ and $h=0$.\\
(i) $J \mapsto m_*(J)$ is strictly increasing on $(\tfrac32,\infty)$.\\
(ii) $J \mapsto \Psi_c(J)$ is strictly decreasing on $(1,\infty)$.\\
(iii) $J\mapsto \Lambda_{t,0}(J)$ is left-continuous at $J=\tfrac32$ for
all $t>\Psi_c(\tfrac32)$. \\
(iv) $J\mapsto \Lambda_{\Psi_c(\bar{J}),0}(J)$ is right-continuous at
$J=\bar{J}$ for all $\bar{J}>\tfrac32$. \\
(v) For every $J\le 3/2$ the map $t \mapsto \Lambda_{t,0}(J)$ is continuous. \\
(vi) For every $J>3/2$ the map $t \mapsto \Lambda_{t,0}(J)$ is continuous except at $t=\Psi_c$ where it exhibits a right-continuous jump.  
\end{theorem}


\subsubsection{Optimal trajectories for $\alpha \in [-1,+1]$, $h \in \R$}
\label{S1.5.4}

\begin{figure}[htbp]
\begin{center}
\begin{tabular}{ccc}
\includegraphics[width=0.3\textwidth]{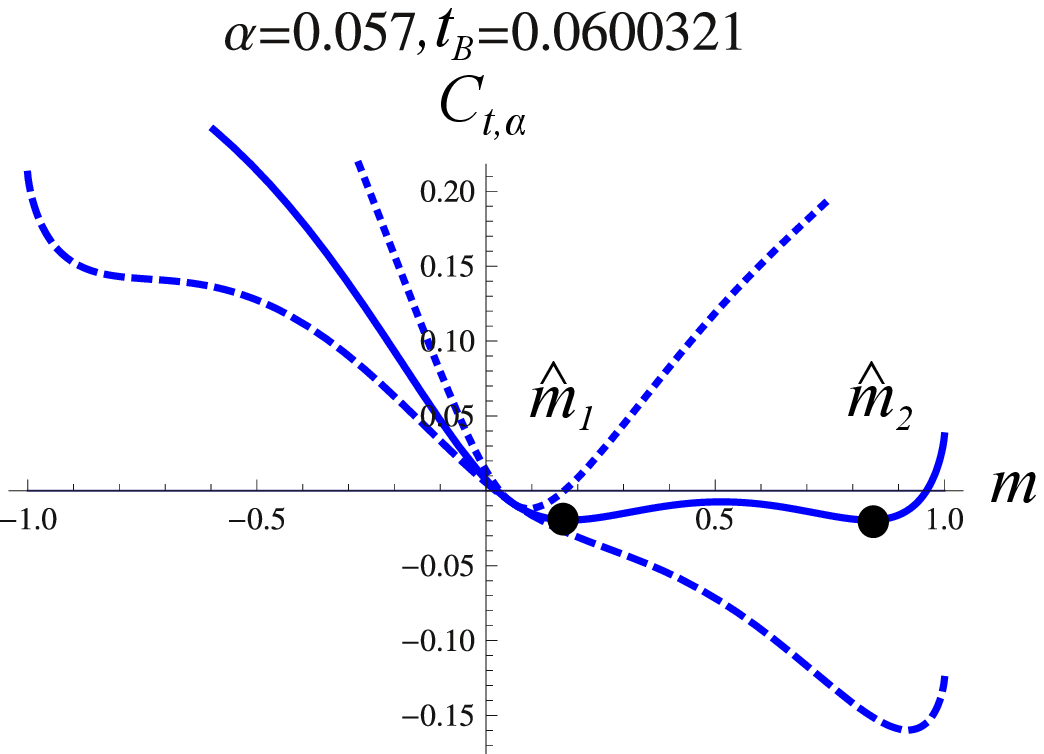} &
\includegraphics[width=0.3\textwidth]{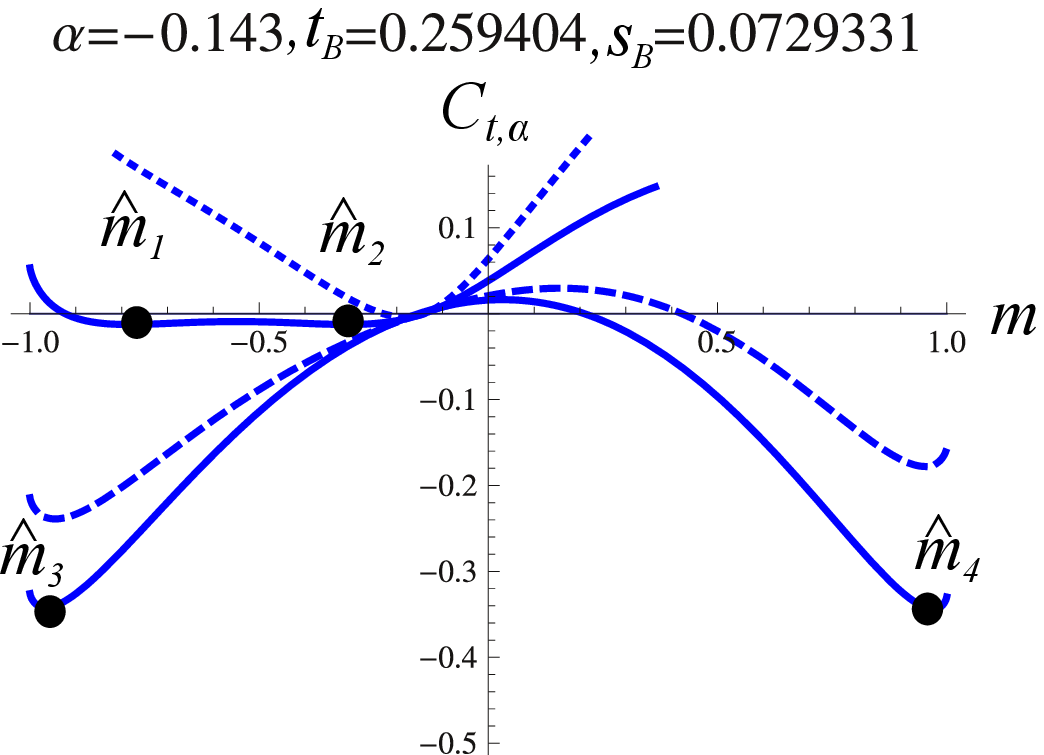} &
\includegraphics[width=0.3\textwidth]{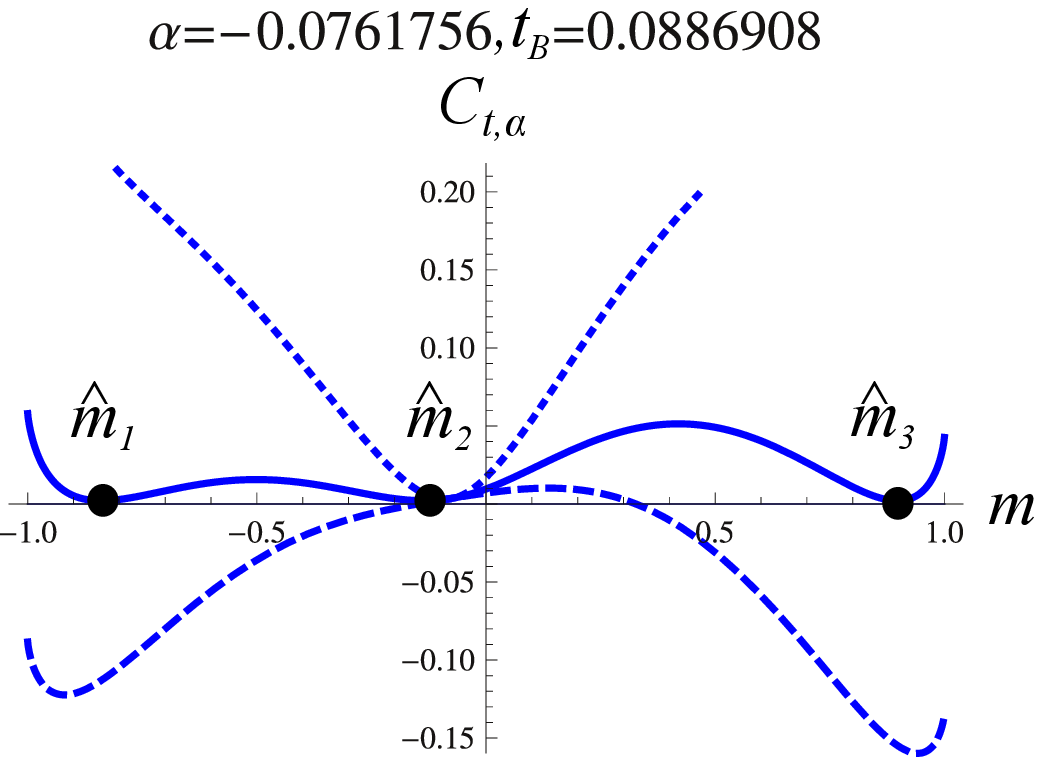}\\
Bifurcation & Double bifurcation & Trifurcation
\end{tabular}
\end{center}
\caption{Different scenarios for the evolution in time of $m \mapsto C_{t,\alpha}(m)$.
Drawn lines: $t=t_B$, $t=t_B,s_B$, $t=t_T$ (times at which multiple global minima 
occur or, equivalently, discontinuity points of $t \mapsto \hat{m}(t,\alpha)$).
Dotted lines: earlier time. Dashed lines: later time.}
\label{fig-3}
\end{figure}

Fore fixed $(J,h)$ and $\alpha$ we say that there is (See Fig.~\ref{fig-3}):
\begin{itemize}
\item 
\emph{No bifurcation} if $\Delta_{t,\alpha} = \{\hat{m}(t,\alpha)\}$, 
for all $t \geq 0$ and the map $t\mapsto \hat{m}(t,\alpha)$ is continuous on $[0,\infty)$.
\item 
\emph{Bifurcation} when there exists a $0<t_B<\infty$ such that
$t \mapsto \hat{m}(t,\alpha)$
continuous except at $t=t_B$ and
$\card{\Delta_{t_B,\alpha}} =2$.
\item 
\emph{Double bifurcation} if there exist times $0<s_B<t_B< \infty$ such that 
$t \mapsto \hat{m}(t,\alpha)$ continuous except 
at $t=t_B$ and $t=s_B$, and
$\card{\Delta_{s_B,\alpha}}=\card{\Delta_{t_B,\alpha}}=2$.
\item \emph{Trifurcation} if there exists a $0<t_T<\infty$ such that
$t \mapsto 
\hat{m}(t,\alpha)$ is continuous except at $t=t_T$ and
$\card{\Delta_{t_T,\alpha} }=3$.
\end{itemize}
The bifurcation times $t_B$ and $s_B$, the trifurcation time $t_T$ and 
the trifurcation magnetization $M_T$ (defined below) all depend on $J,h$.

The following theorem summarizes the behaviour of $\Delta_{t,\alpha}$ (and 
therefore of the minimizing trajectories $\hat{\phi}_{t,\alpha}$) for different 
$t,\alpha$. For $J>\tfrac32$, let 
\begin{equation}
\label{F}
\begin{aligned}
F(m) &:= \dfrac{mk'^{J,h}(m)-k^{J,h}(m)}{\csch[\acoth(k'^{J,h}(m))]},\\
U_B=U_B(J,h) &:=\max\limits_{m \in [0,1]} F(m),\\
L_B=L_B(J,h) &:=\min\limits_{m \in [-1,0]} F(m).
\end{aligned}
\ee

\begin{theorem}
\label{generaltheorem}
{\rm (See Figs.~\ref{fig-3}--\ref{fig-2}.)}
\begin{itemize}
\item[{\rm (1)}]
Suppose that $k^{J,h}(\alpha) \neq 0$.\\
{\rm (1a)} 
If $k^{J,h}(\alpha)>0$ and $\alpha>0$, then there are $m_R^+>0$ and $t_R=t_R(m_R^+)>0$ 
(implicitly calculable from \eqref{mRplus}) such that $t\mapsto\hat{m}(t)$ is 
strictly increasing on $[0,t_R]$ and strictly decreasing on $[t_R,\infty)$ with 
$\hat{m}(t_R)=m_R^+>m^\infty$.\\
{\rm (1b)} 
If $k^{J,h}(\alpha)<0$ and $\alpha>0$, then $t\mapsto\hat{m}(t)$ is strictly decreasing 
on $[0,\infty)$.\\
{\rm (1c)} 
If $k^{J,h}(\alpha)>0$ and $\alpha<0$, then $t\mapsto\hat{m}(t)$ is strictly increasing 
on $[0,\infty)$.\\
{\rm (1d)} 
If $k^{J,h}(\alpha)<0$ and $\alpha<0$, then there are $m_R^->0$ and $t_R=t_R(m_R^-)>0$ 
(implicitly calculable from \eqref{mRminus}) such that $t\mapsto\hat{m}(t)$ is 
strictly decreasing on $[0,t_R]$ and strictly increasing on $[t_R,\infty)$ with 
$\hat{m}(t_R)=m_R^-<\alpha$.\\
In all cases $\hat{m}(0)=\alpha$ and $\lim_{t\to\infty} \hat{m}(t)=m^\infty$. 
\item[{\rm (2)}]
Suppose that $h=0$.\\
{\rm (2a)} 
If $0<J\leq 1$, then there is no bifurcation.\\
{\rm (2b)} 
If $1<J \leq \tfrac32$, then there is bifurcation only for $\alpha=0$.\\
{\rm (2c)}
If $J>\tfrac32$, then there is bifurcation if $\alpha \in (-U_B,U_B)$ and no bifurcation otherwise.
\item[{\rm (3)}] 
Suppose that $h>0$.\\
{\rm (3a)} 
If $0<J\leq 1$, then there is no bifurcation.\\
{\rm (3b)} 
If $1<J\leq \tfrac32$, then there is bifurcation for $\alpha \in [-1,U_B)$ and no 
bifurcation for $\alpha \in [U_B,1]$.\\
{\rm (3c)} 
If $J>\tfrac32$, then there exists a $h_*=h_*(J)>0$ such that 
\begin{list}{-}{\addtolength{\itemsep}{-1mm} \addtolength{\topsep}{-2mm}}
\item for every $0< h<h_*$ there exists a $M_T \in (L_B,U_B)$ with $M_T<0$ such that 
there is
\begin{list}{*}{\addtolength{\itemsep}{-1mm} \addtolength{\topsep}{-2mm}}
\item no bifurcation for $\alpha \in [U_B,1]$,
\item bifurcation for $\alpha \in (M_T,U_B)$,
\item trifurcation for $\alpha=M_T$,
\item double bifurcation for $\alpha \in (L_B,M_T)$,
\item bifurcation for $\alpha \in [-1,L_B]$.
\end{list}
\item for every $h \geq h^*$ the behavior is the same as in {\rm (3b)}.
\end{list}
In all cases $\alpha \mapsto t_B(\alpha)$ is continuous and decreasing and 
$\alpha \mapsto s_B(\alpha)$ is continuous and increasing.
\end{itemize}
\end{theorem}

\noindent 
Theorem~\ref{generaltheorem} gives a complete picture of the bifurcation scenario. 
Regime (1) ---which includes cases with zero and nonzero magnetic field--- describes two types of behavior of optimal magnetization trajectories:  monotone trajectories [cases (1b) and (1c)] and trajectories with \emph{overshoot} [cases (1a) and (1d)].  In the latter, $\hat{m}(t)$ increases to some magnetization 
$m_R^+$ larger ($m_R^-$ smaller) than $m^\infty$ and afterwards decreases (increases) 
to $m^\infty$.  
Regimes (2)and (3) refer to the existence of bifurcations and trifurcations.  We observe that the 
different bifurcation behaviors ---no bifurcation, single and double bifurcation--- hold for whole intervals of the conditioning magnetization.  In contrast, trifurcation appears at a single final magnetization for each $h\neq 0$.


\begin{figure}[hbtp]
\label{Overshoot}
\begin{center}
\begin{tabular}{cc}
\includegraphics[width=0.45\textwidth]{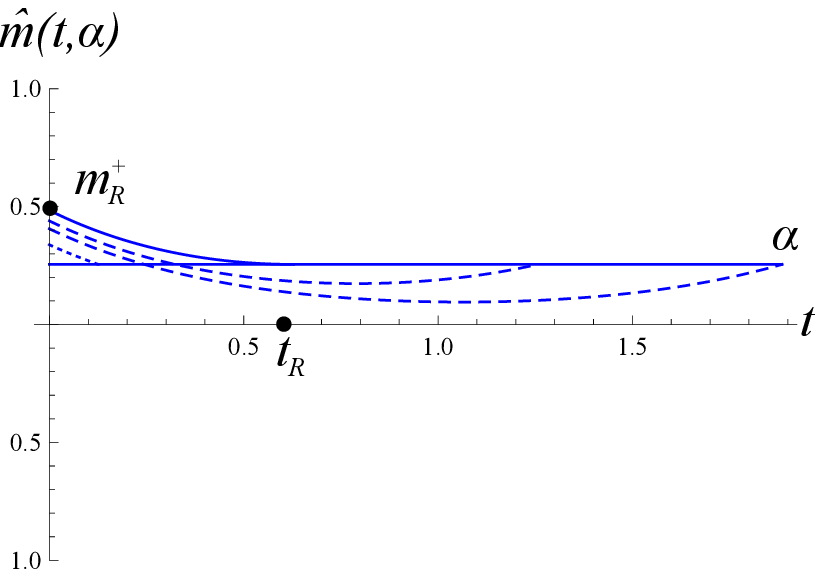} &
\includegraphics[width=0.45\textwidth]{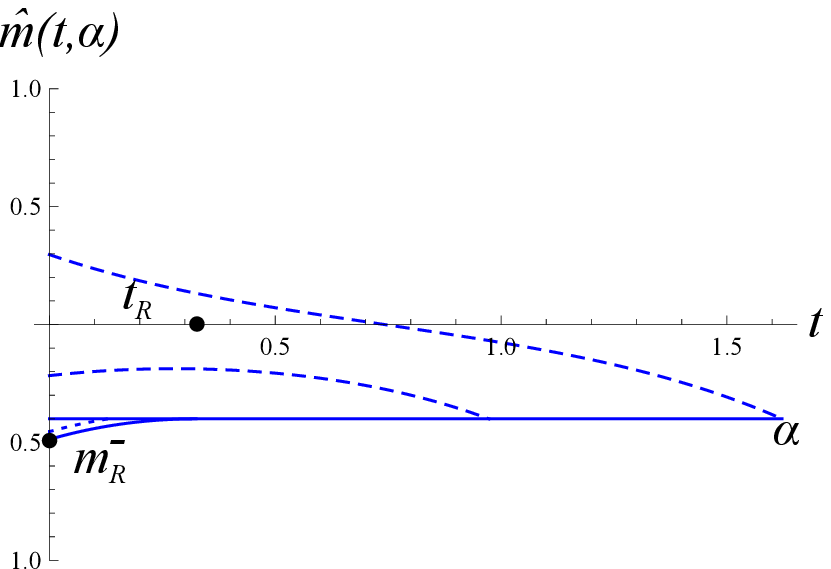} \\
Regime (1a) & Regime (1d)\\
\includegraphics[width=0.45\textwidth]{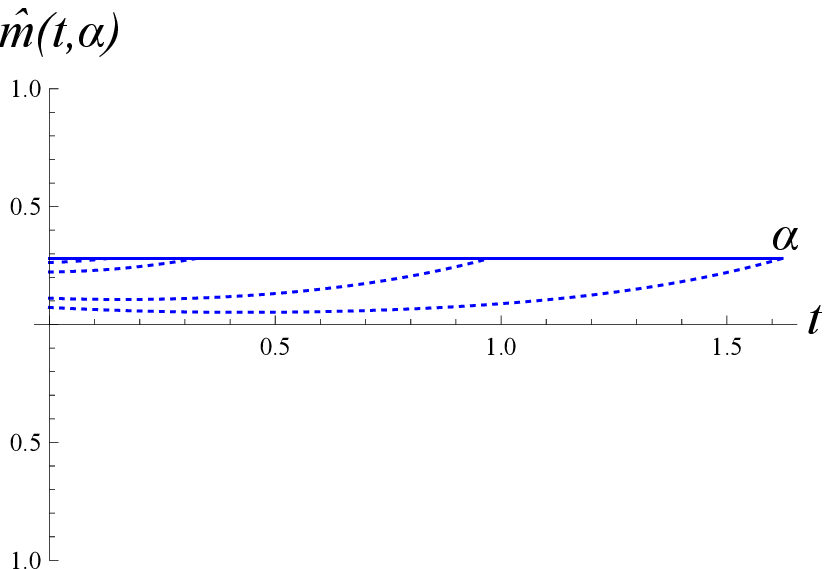} &
\includegraphics[width=0.45\textwidth]{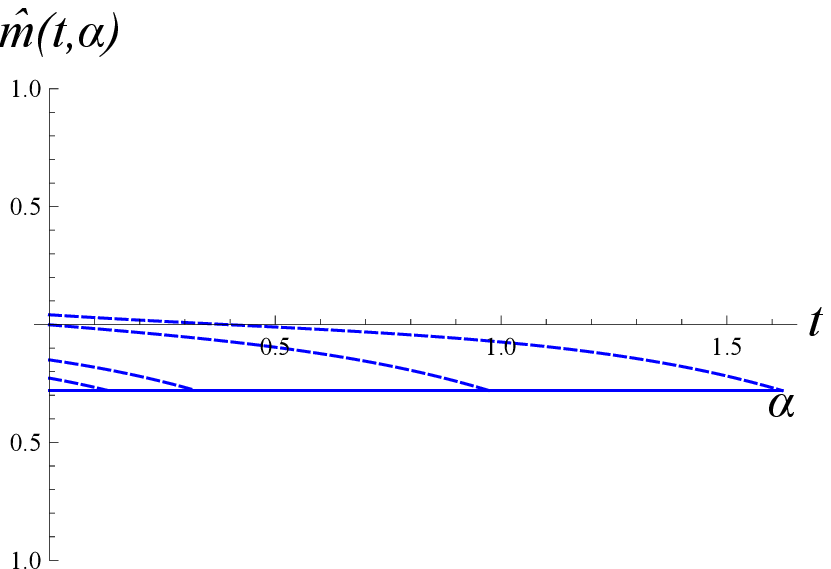}\\
Regime (1b) & Regime (1c)
\end{tabular}
\end{center}
\caption{Different regimes of Theorem \ref{generaltheorem}. Evolution in time 
of the minimizing trajectories $\pm(\hat{\phi}_{t,\alpha}(s))_{0\leq s \leq t}$ 
for $t<t_R$ (dotted), $t=t_R$ (drawn), $t>t_R$ (dashed).}
\label{fig-2}
\end{figure}

\subsubsection{Gibbs versus non-Gibbs}
\label{S1.5.5}

Theorem~\ref{nongibbs-bifurcation} establishes the equivalence of bifurcation and 
discontinuity of specifications, as proposed in the program put forward in 
\cite{vEFedHoRe10}.  Due to this equivalence, the following corollary provides a full characterization of the different Gibbs--nonGibbs scenarios appearing during the infinite-temperature evolution of the Curie-Weiss model. Let
\be
0< \Psi_U:=t_B(U_B) < \Psi_T:=t_B(M_T) 
< \Psi_L:=t_B(L_B) < \Psi_*:=t_B(-1),
\ee
and let $M_B$ be the solution of $t_B(M_B)=\Psi_L$. Denote $\cD_t \subset [-1,+1]$ 
the set of $\alpha$-values for which $\alpha \mapsto \gamma_t(\cdot|\alpha)$ is
discontinuous. 

\begin{corollary}
\label{GNG-coro}
{\rm (See Fig.~\ref{PD}.)}
\begin{itemize}
\item[{\rm (1)}] 
Let $h=0$.\\
{\rm (1a)} If $0<J\leq 1$, the evolved measure $\mu_t$ is Gibbs for all $t \geq 0$.\\
{\rm (1b)} If $1<J\leq \tfrac32$, then $\mu_t$ is 
\begin{list}{-}{\addtolength{\itemsep}{-1mm} \addtolength{\topsep}{-2mm}}
\item Gibbs for $0 \leq t \leq \Psi_c$,
\item non-Gibbs for $t > \Psi_c$ with $\cD_t=\{0\}$.
\end{list}
{\rm (1c)} If $J > \tfrac32$, then $\mu_t$ is 
\begin{list}{-}{\addtolength{\itemsep}{-1mm} \addtolength{\topsep}{-2mm}}
\item Gibbs for $0 \leq t \leq \Psi_U$,
\item non-Gibbs for $t > \Psi_U$ with 
\begin{list}{*}{\addtolength{\itemsep}{-1mm} \addtolength{\topsep}{-2mm}}
\item $\cD_t=\{\pm \alpha\}$ for some $\alpha \in (-U_B,U_B)$ if $\Psi_U < t < \Psi_c$,
\item $\cD_t=\{0\}$ if $t \geq \Psi_c$.
\end{list}
\end{list}
\item[{\rm (2)}] 
Let $h>0$.\\
{\rm (2a)} If $0<J\leq 1$, then $\mu_t$ is Gibbs for $t \geq 0$.\\
{\rm (2b)} If $1<J\leq \tfrac32$, then $\mu_t$ is 
\begin{list}{-}{\addtolength{\itemsep}{-1mm} \addtolength{\topsep}{-2mm}}
\item Gibbs for $0 \leq t \leq \Psi_U$,
\item non-Gibbs for $\Psi_U<t \leq \Psi_*$ with $\cD_t=\{\alpha\}$ for 
some $\alpha \in [-1,U_B)$,
\item Gibbs for $t> \Psi_*$.
\end{list}
{\rm (2c)} If $J>\tfrac32$ and $h<h^*$ small enough, then $\mu_t$ is 
\begin{list}{-}{\addtolength{\itemsep}{-1mm} \addtolength{\topsep}{-2mm}}
\item Gibbs for $0 \leq t \leq \Psi_U$,
\item non-Gibbs for $\Psi_U<t \leq \Psi_*$ with 
\begin{list}{*}{\addtolength{\itemsep}{-1mm} \addtolength{\topsep}{-2mm}}
\item $\cD_t=\{\alpha\}$ for some $\alpha \in [M_B,U_B)$ if $\Psi_U<t \leq \Psi_L$,
\item $\cD_t=\{\alpha^1,\alpha^2\}$ for some $\alpha^1,\alpha^2 \in 
(L_B,M_B)$ if $\Psi_L < t < \Psi_T$,
\item $\cD_t=\{\alpha\}$ for some $\alpha \in [-1,M_T]$ if 
$\Psi_T \leq t \leq \Psi_*$.
\end{list}
\item Gibbs for $t> \Psi_*$.
\noindent
If $h \geq h^*$, then the behaviour is as in {\rm (2b)}.
\end{list}
\end{itemize}
In all cases $\alpha^1,\alpha^2,\alpha$  depend on $(t,J,h)$.
\end{corollary}

\begin{figure}[h]
\begin{center}
\begin{tabular}{cc}
\includegraphics[width=0.45\textwidth]{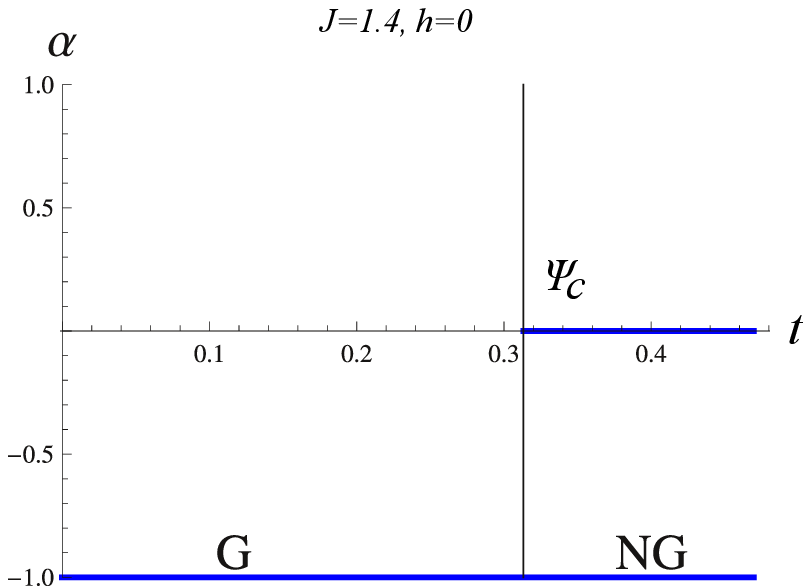} &
\includegraphics[width=0.45\textwidth]{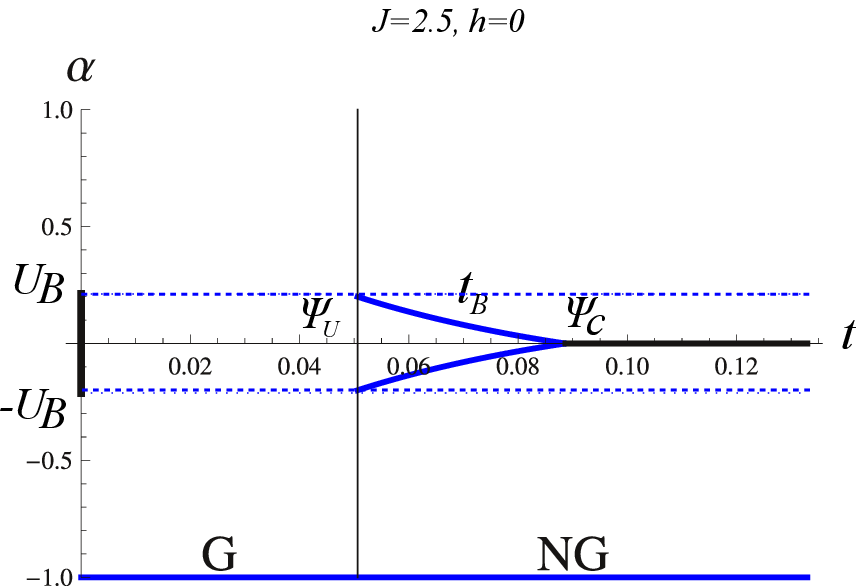}\\
\includegraphics[width=0.45\textwidth]{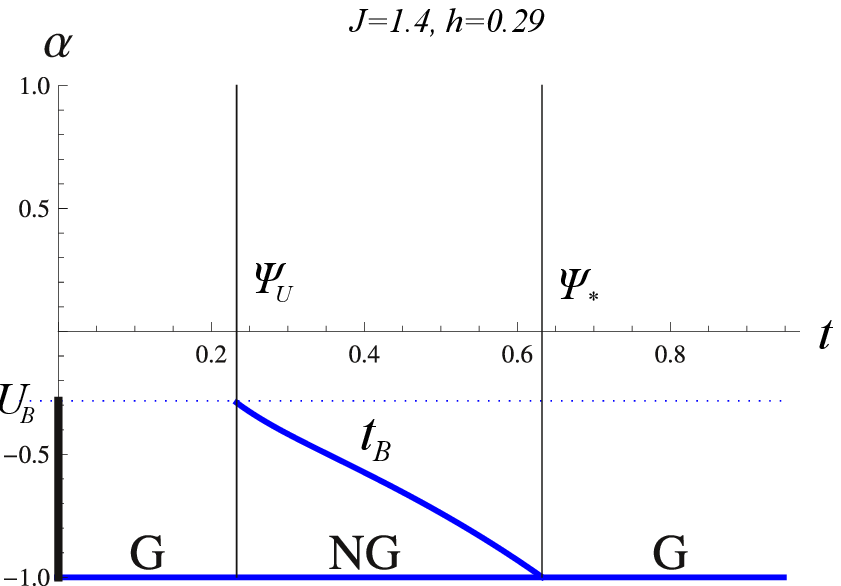} & 
\includegraphics[width=0.45\textwidth]{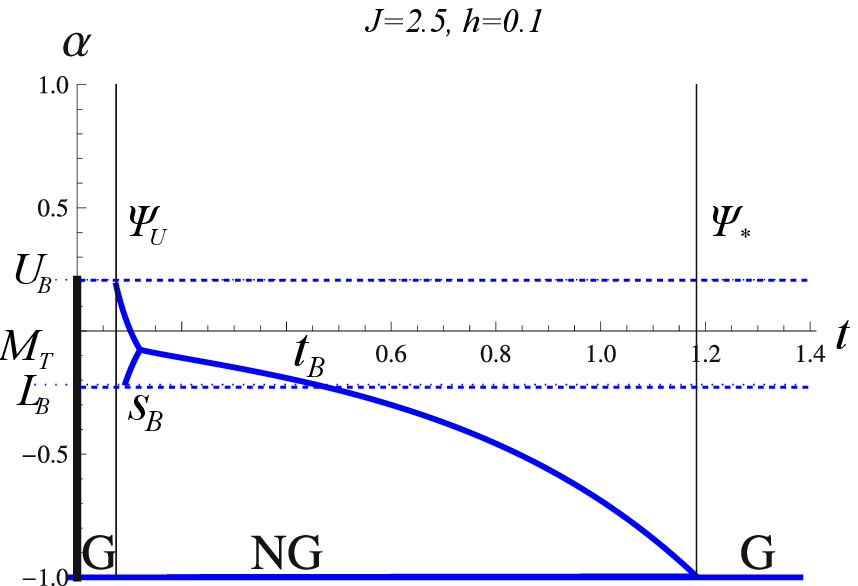}\\
&\\
{$1<J \leq \tfrac{3}{2}$} & { $J>\tfrac{3}{2}$}
\end{tabular}
\end{center}
\caption{Summary of Corollary \ref{GNG-coro}: Time versus bad magnetizations for 
different regimes. On the  vertical $\alpha$-axis, indicated by a thick line, is
the set of bad magnetizations. G=Gibbs, NG=non-Gibbs.}
\label{PD}
\end{figure}


\section{Proof of Proposition~\ref{Thereduction} and 
Theorems~\ref{nongibbs-bifurcation}--\ref{Cone-Monotonicity}}
\label{S2}

Proposition~\ref{Thereduction} is proven in Section~\ref{S2.1}, 
Theorems~\ref{nongibbs-bifurcation}--\ref{Cone-Monotonicity} are proven in
Sections~\ref{S2.2}--\ref{S2.4}.


\subsection{Proof of Proposition~\ref{Thereduction}}
\label{S2.1}

\bpr
First note that, by \eqref{originalcost},   
\be
\label{twoinf}
\inf_{\substack{\phi\colon\,\phi(t)=\alpha}} I^t(\phi)
= \inf\limits_{m \in [-1,+1]} \left\{I_S(m) 
+ \inf_{\substack{\phi\colon\,\phi(0)=m, \\ \phi(t)=\alpha}} I^t_D(\phi)\right\}.
\ee
It follows from (\ref{originalcost}--\ref{actionid}) and the calculus of variations 
that the stationary points of the right-hand side of \eqref{twoinf} are given by the 
Euler-Lagrange equation, complemented with a free-left-end condition and a 
fixed-right-end condition:
\be
\label{system}
\begin{split}
\frac{\partial}{\partial s} \frac{\partial L}{\partial \dot m}(\phi(s),\dot{\phi}(s)) 
&= \frac{\partial L}{\partial m}(\phi(s),\dot{\phi}(s)), \qquad s \in \ (0,t),\\
\pd{L}{\dot m}(\phi(s),\dot{\phi}(s))\Big|_{s=0} 
&= \pd{I_S}{m}(\phi(s))\Big|_{s=0},\\
\phi(t)& =\alpha.
\end{split}
\ee
The first and the third equation in \eqref{system} come from the third infimum in 
\eqref{twoinf} and, together with \eqref{Lform}, determine the form 
\eqref{trajectory} of the stationary trajectory.
Inserting this form into \eqref{originalcost} we identify
\begin{equation}
\label{Cbhf}
I^t(\hat{\phi}^m_{t,\alpha})  \;=\; C_{t,\alpha}(m),
\end{equation}
as stated in \eqref{cost}--\eqref{cost.1}. 
This identity reduces \eqref{cvp} to a 
one-dimensional variational problem,
\be
\label{1dp}
\inf_{\substack{\phi\colon\,\phi(t)=\alpha}} I^t(\phi)
= \inf_{m \in [-1,+1]} I^t(\hat{\phi}^m_{t,\alpha}) 
= \inf_{m \in [-1,+1]} C_{t,\alpha}(m)
\ee 
The second equation in \eqref{system} 
corresponds to the second infimum in \eqref{twoinf} or, equivalently, to the rightmost infima in \eqref{1dp}. It gives a trade-off between the static 
and the dynamic cost, establishing a relation between the initial magnetization and 
the initial derivative. After some manipulations this equation can be written in the 
form
\be
\label{iequation}
-\tfrac12\,q = a^J(m)\cosh(2h)+ b^J(m)\sinh(2h), \qquad m = \hat{\phi}^m_{t,\alpha}(0),
\quad q = \dot{\hat{\phi}}^m_{t,\alpha}(0).
\ee
Differentiating \eqref{trajectory}, we get 
\be
\dot{\hat{\phi}}^m_{t,\alpha}(s) 
= 2 \csch(2t) \Big\{\alpha \cosh(2s) - m \cosh(2(t-s))\Big\},
\ee
and eliminating $q$ from the last identity in \eqref{iequation} in favor of $t$ and $\alpha$, we conclude that $m$ must be a solution of
\eqref{Pequation}. 
Imposing this restriction to the chain of identities \eqref{1dp} we obtain \eqref{simplification}
and hence \eqref{trajectory}--\eqref{trajectory.00}.
\epr

From now on our arguments rely on the study of \eqref{Pequation}, combined with
continuity properties of $C_{t,\alpha}$ as a function of $t,\alpha$.


\subsection{Proof of Theorem~\ref{nongibbs-bifurcation}}
\label{S2.2}

Equation \eqref{minvsspecif} follows in the same way as in the proof of Theorem 2.5 in 
\cite{ErKu10}. Having disposed of this identity, we can now proceed to prove the 
equivalence. The proof relies on the following lemma.

\begin{lemma}
For any $t>0$ and $\alpha_0 \in [-1,+1]$, there exists an open neighbourhood $\mathcal{N}_{\alpha_0}\neq \emptyset$ of the later, such that for all $\alpha \in \mathcal{N}_{\alpha_0}\setminus \{\alpha_0\}$
\begin{enumerate}
\item $C_{t,\alpha}$ has only one global minimum, namely, $\hat{m}(t,\alpha)$.
\item $\tilde{\alpha}\mapsto \hat{m}(t,\tilde{\alpha})$ is continuous at $\alpha$. If $C_{t,\alpha_0}$ has a unique global minimum, the continuity is also valid at $\alpha=\alpha_0$.
\item If $C_{t,\alpha_0}$ has multiple global minima, for two of them, namely, $\hat{m}_A(t,\alpha_0)$ and $\hat{m}_B(t,\alpha_0)$
$$\lim\limits_{\alpha \downarrow \alpha_0} \hat{m}(t,\alpha)=\hat{m}_A(t,\alpha_0)
\quad \text{and} \quad 
\lim\limits_{\alpha \uparrow \alpha_0} \hat{m}(t,\alpha)=\hat{m}_B(t,\alpha_0).
$$
\end{enumerate}
\end{lemma}
\bpr
A straightforward study of \eqref{Pequation} shows that $C_{t,\alpha}$ has a finite number of critical points, for every fixed choice of $J,h,t,\alpha$.

Clearly, $\alpha \mapsto C_{t,\alpha}$ and $\alpha \mapsto l_{t,\alpha}$ are continuous 
with respect to the infinity norm in $C([-1,+1],\mathbb{R})$. This, together with the 
fact that the left-hand side of \eqref{Pequation} does not depend on $\alpha$, implies continuity of any critical point with respect to $\alpha$.

Let $\hat{m}_i(t, \alpha_0), \ i=1,\ldots,v$ be the global minima of $C_{t,\alpha_0}$.
By continuity of the critical points, there exists a neighbourhood $\widetilde{\mathcal{N}}_{\alpha_0}$ and smooth functions $\widetilde{\mathcal{N}}_{\alpha_0} \ni \alpha \mapsto \overline{m}_i(t,\alpha)$,  $i=1,\ldots,v$, such that 
\begin{itemize}[topsep=0pt, partopsep=0pt,itemsep=0pt]
\item[i.] $\overline{m}_i(t, \alpha)$ are local minima of $C_{t,\alpha}$,
\item[ii.] $\lim\limits_{\alpha \to \alpha_0}\overline{m}_i(t, \alpha)= \hat{m}_i(t, \alpha_0)$.
\end{itemize}
This properties proves the lemma if $v=1$. Otherwise, let
$$B_i(\alpha):= C_{t,\alpha}(\overline{m}_i(t,\alpha)).$$
The minimal cost is attained at the smallest of them:
 \[C_{t, \alpha}(\hat{m}(t, \alpha))=\min_i B_i(\alpha)\;.\]
Note that there is coincidence at $\alpha_0$ due to the assumed multiplicity of minima:
\be
B(\alpha_0) \;:=\: B_i(\alpha_0)\quad,\quad i=1,\ldots v\;.
\ee
\label{eq:bb0}
We expand the functions $B_i$ up to first order order 
\be
\label{eq:bb1}
B_i(\alpha)= B(\alpha_0) + B_i'(\alpha_0) (\alpha- \alpha_0) + O(\alpha- \alpha_0)\;,
\ee
and observe that,  
\be
\label{eq:bb2}
B_i'(\alpha_0) \neq B_j'(\alpha_0), \qquad i \neq j.
\ee
The latter is due to the strict monotonicity of $\pd{C_{t,\alpha}}{\alpha}$ and the fact that each $\overline{m}_i(t,\alpha)$ is a critical points of the function $C_{t,\alpha}(\,\cdot\,)$.
From \eqref{eq:bb0}--\eqref{eq:bb2} we conclude that for $\alpha$ in a possibly smaller neighbourhood $\mathcal{N}_{\alpha_0}\subset\widetilde{\mathcal{N}}_{\alpha_0}$ there is a unique global minimum, and
that property $3.$ holds with
$$a=\arg\min_i B_i'(\alpha_0), \qquad b=\arg\max_i B_i'(\alpha_0),$$
and 
$$\hat{m}_A(t,\alpha_0):=\hat{m}_a(t,\alpha_0), \qquad \hat{m}_B(t,\alpha_0):=\hat{m}_b(t,\alpha_0).$$
\epr

We are ready to prove Theorem~\ref{nongibbs-bifurcation}.
\bpr
Suppose that $C_{t,\alpha_0}$ has a unique minimizer, denoted by 
$\hat{m}(t,\alpha_0)$ and let $\mathcal{N}_{\alpha_0}$ be the neighbourhood of the previous lemma. 
Then \eqref{spec-initialm} holds for every $\alpha \in \cN_{\alpha_0}$, and the continuity of $m \mapsto \Gamma_t(z,m)$ for every $t,z$ gives the desired continuity of $\alpha \mapsto \gamma_t(\cdot \mid \alpha)$ at $\alpha=\alpha_0$.

To prove necessity, assume that $C_{t,\alpha_0}$ has multiple global minima. Consider $\hat{m}_A$ and $\hat{m}_B$ as in the previous lemma. Then, we have that there exist sequences $\alpha^-_n<\alpha_0< \alpha^+_n$ converging to $\alpha_0$ and such that 
$\gamma_t(\cdot \mid \alpha_n^{\pm})=\Gamma_t(\cdot,\hat{m}(t, \alpha_n^\pm))$ and
\be
\lim_{n \to \infty}\hat{m}(t,\alpha_n^-)=\hat{m}_B(t,\alpha_0) 
\neq \hat{m}_A(t,\alpha_0) = \lim_{n \to \infty}\hat{m}(t,\alpha_n^+). 
\ee
Again using continuity of $\Gamma_t$ with respect to $m$, we get
$$\lim\limits_{n\to \infty} \gamma_t(z \mid \alpha^-_n)=\Gamma_t(z,\hat{m}_B(t,\alpha_0))\neq
\Gamma_t(z,\hat{m}_A(t,\alpha_0))=\lim\limits_{n\to \infty} \gamma_t(z \mid \alpha^+_n).$$
Hence $\alpha_0$ is a bad magnetization.
\epr

\subsection{Proof of Theorem~\ref{hmfzero}}
\label{S2.3}

To determine which solutions of \eqref{Pequation} are global minima of $C_{t,0}$ when 
$h=0$, we will pursue the following strategy. Using \eqref{Pequation} we can 
write $t$ as a function of $m$:
\be
\label{tm}
t(m) := \tfrac12 \acoth \left(\frac{a^J(m)}{m}\right).
\ee
This allows us to determine for which time $t$ the magnetization $m$ can be a \emph{possible 
minimum} (i.e., a solution of \eqref{Pequation}). 

\begin{lemma}
\label{COSTOM}
Let $A\subseteq [-1,+1]$ be the set of $m$-values such that $m$ is the solution of 
\eqref{Pequation} for some $t>0$, i.e., $A = \{m\in[-1,+1]\colon\,a^J(m)/m>1\}$. Then, 
for every $m \in A$,
\be
\label{Mcosto}
C_{t(m),0}(m) = \tfrac12 J m^2 + \tfrac12 \log\,[1-m \tanh(Jm)] =: C_M(m).
\ee
In words, \eqref{Mcosto} is the cost for $m$ at the time at which it is a possible 
minimum.
\end{lemma}

\bpr
Insert \eqref{tm} into \eqref{cost} and use \eqref{system}.
\epr

We now start the proof of Theorem~\ref{hmfzero}. 

\begin{proof}
First note that $l_{t,0}$ is linear with slope $\coth(2t) \in (1,\infty)$ and 
$l_{t,0}(0)=0$, and that $a^J$ is antisymmetric. Hence, if $m$ is a solution 
of \eqref{Pequation}, then also $-m$ is a solution. Further note that  
\be
\begin{split}
{a^J}'(m) &= (2J-1) \cosh(2Jm) - 2Jm \sinh(2Jm),\\
{a^J}'' (m)&= 4J(J-1) \sinh(2Jm) - 4J^2m \cosh(2Jm).
\end{split}
\ee

\medskip\noindent
(i)
If $0<J<\tfrac12$, then ${a^J}'(m)<0$ for all $m$, and hence $m=0$ is the 
unique solution for all $t>0$. If $\tfrac12 \leq J \leq 1$, then ${a^J}''(m)<0$ 
for all $m$, hence $a^J$ is convex, and so it suffices to compare slopes at 0: 
${k^{J,0}}'(0)={a^J}'(0)=2J-1<1$ and ${l_{t,0}}'(0)=\coth(2t)>1$. Again, $m=0$ 
is the unique solution for all $t>0$ (see Fig.~\ref{gb1}).

\medskip\noindent
(ii)	
As before, ${a^J}''(m)<0$ for all $m$, but now the slopes at $0$ can be equal, 
which occurs when $t = \Psi_c$ with $\Psi_c$ defined in \eqref{CRITICALT}. This proves 
that $\Delta_t = \{0\}$ for $0<t\leq \Psi_c$ and $\Delta_t \subseteq \{-\hat{m}(t),0,
\hat{m}(t)\}=$ the set of solutions of \eqref{Pequation} for $t>\Psi_c$. It is easily 
seen from Fig.~\ref{gb2} that $\hat{m}(t)$ is continuous and strictly increasing 
on $[\Psi_c,\infty)$ and $\hat{m}(\Psi_c)=0$. It remains to show that $\{-\hat{m}(t),
\hat{m}(t)\}$ are the global minima for all $t>\Psi_c$. This follows from the strategy 
behind the proof of Lemma~\ref{COSTOM}. Since $C_{t,0}(0)=0$ for all $t>0$, it 
suffices to prove that $m \mapsto C_M(m)$ is strictly decreasing. From \eqref{Mcosto} 
we have
\be
\label{Dcm}
C_M'(m) = \frac{\partial C_{t(m),0}}{\partial m}(m)
+ \frac{\partial C_{t(m),0}}{\partial t}(m)\,t'(m).
\ee
The first term is zero by the definition of $t(m)$ (each $m$ is a stationary point 
of $C_{t,0}$ at time $t=t(m)$). The second term is $<0$ because $t'(m)>0$ and 
\be
\label{dC/dt}
\begin{aligned}
\frac{\partial C_{t,0}}{\partial t}(m)
&= L(\hat{\phi}^m_{t,0}(t),\dot{\hat{\phi}}^m_{t,0}(t))\\ 
&\qquad + \int\limits_0^t \left[\frac{\partial L}{\partial m}
\big(\hat{\phi}^m_{t,0}(s),\dot{\hat{\phi}}^m_{t,0}(s)\big)\,
\frac{\partial\hat{\phi}^m_{t,0}}{\partial t}(s) 
+ \frac{\partial L}{\partial \dot m}\big(\hat{\phi}^m_{t,0}(s),
\dot{\hat{\phi}}^m_{t,0}(s)\big)\,
\frac{\partial \dot{\hat{\phi}}^m_{t,0}}{\partial t}(s)\right] ds\\
&= L(0,\dot{\hat{\phi}}^m_{t,0}(t)) 
+ \int\limits_0^t \frac{\partial}{\partial s}
\left\{\frac{\partial L}{\partial \dot m}
\big(\hat{\phi}^m_{t,0}(s),\dot{\hat{\phi}}^m_{t,0}(s)\big)\, 
\frac{\partial \hat{\phi}^m_{t,0}}{\partial t}(s)\right\}ds\\ 
&= L(0,r) - \pd{L}{\dot m}(0,r)\,r\\
&= - \tfrac12 \sqrt{4 + r^2} + 1 < 0
\end{aligned}
\ee
with $r=\dot{\hat{\phi}}^m_{t,0}(t)$, where the second equality uses \eqref{system}.
Since $C_M(0)=0$, this yields the claim (see Figs.~\ref{gb2},\ref{Ct2},\ref{CM2}).

\medskip\noindent
(iii)
This case is more difficult, because $a^J$ no longer is convex on $(0,1)$. Let 
$\Psi_c^1$ be the first time at which a solution different from $0$ exists. To 
identify $\Psi_c^1$, let
\be
T_{m}(x):=(x-m){a^J}'(m) + a^J(m),
\ee
and let $m_1$ be the solution of the equation $T_{m_1}(0) = 0 =-m_1 {a^J}'(m_1) 
+ a^J(m_1)$, i.e.,
\be
\label{Tc}
m_1 = \frac{a^J(m_1)}{{a^J}'(m_1)}.
\ee
From $m_1$ we get $t^1_c$ by using \eqref{tm}: $t^1_c=t(m_1)$. As before, a solution of 
\eqref{Pequation} for $t \geq t^1_c$ is not necessarily a minimum. To find out when it 
is, we follow the same strategy as in case (ii). Again, $C_{t,0}(0)=0$ for all $t>0$, 
and hence we must look for $m_*>0$ such that 
\begin{equation}
\label{mjump}
C_M(m_*)=0.
\end{equation}
Knowing $m_*$, we are able to compute $\Psi_c$ using \eqref{tm},
\begin{equation}
\label{tjump}
t_*=t(m_*).
\end{equation}
In words, $t_*$ is the first time at which $0$ no longer is a minimum. As in case (ii), 
it suffices to prove that $m\mapsto C_M(m)$ is strictly decreasing on $(m_*,\infty)$.
Again, we have \eqref{Dcm}. Since
\be
t'(m) = \frac12 (\atanh )' \left( \frac{a^J(m)}{m}\right) 
\left\{{a^J}'(m) - \frac{a^J(m)}{m}\right\}\, \frac{1}{m},
\ee
it follows that $t'(m)=0$ if and only of $m= a^J(m)/{a^J}'(m)$, which is the same condition 
as \eqref{Tc}. This gives us a graphical argument to conclude that $t'(m)<0$ for $0<m<m_1$ 
and $t'(m)>0$ for $m>m_1$ (see Figs.~\ref{gb3},\ref{Ct3},\ref{CM3}). On the other hand, 
$m_*>m_1$.
\end{proof}

\begin{figure}[htbp]
\begin{center}
\includegraphics[width=0.4 \textwidth]{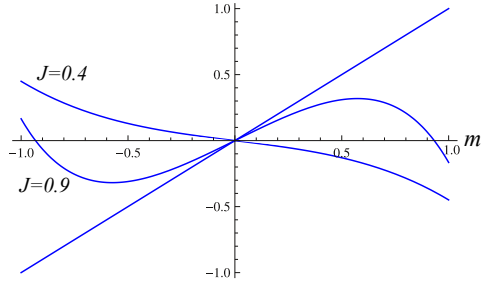}
\end{center}
\caption{$m \mapsto a^J(m)$, $m \mapsto l_{t,0}(m)$ for Regime (i).}
\label{gb1}
\end{figure}

\begin{figure}[htbp]
\begin{center}
\begin{tabular}{ccc}
\subfigure[]{\label{gb2}\includegraphics[width=0.3 \textwidth]{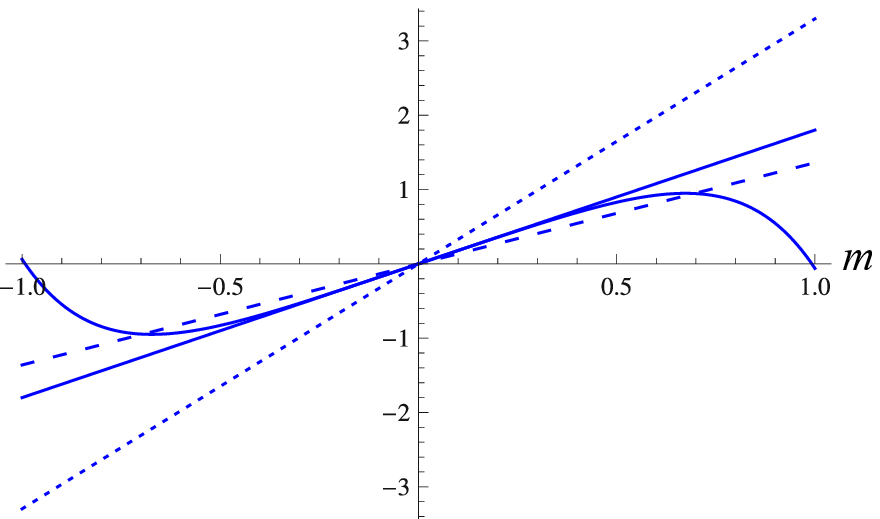}}&
\subfigure[]{\label{Ct2}\includegraphics[width=0.3 \textwidth]{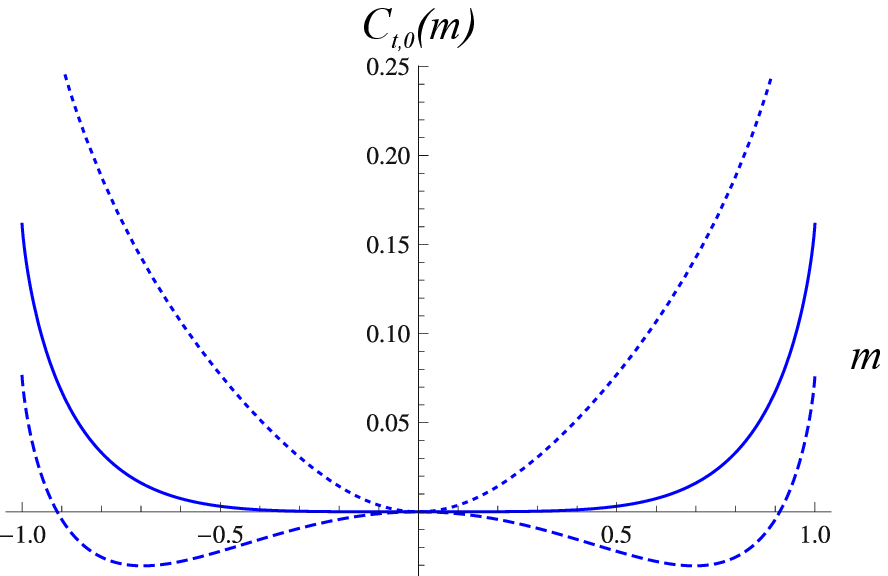}}&
\subfigure[]{\label{CM2}\includegraphics[width=0.3 \textwidth]{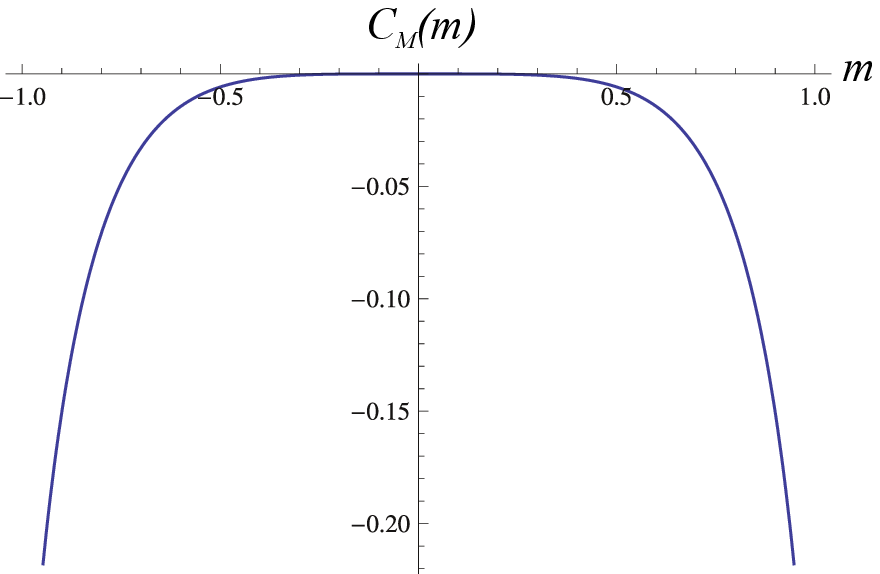}}\\
& Regime (ii), $J=1.3$. & \\
\subfigure[]{\label{gb3}\includegraphics[width=0.3 \textwidth]{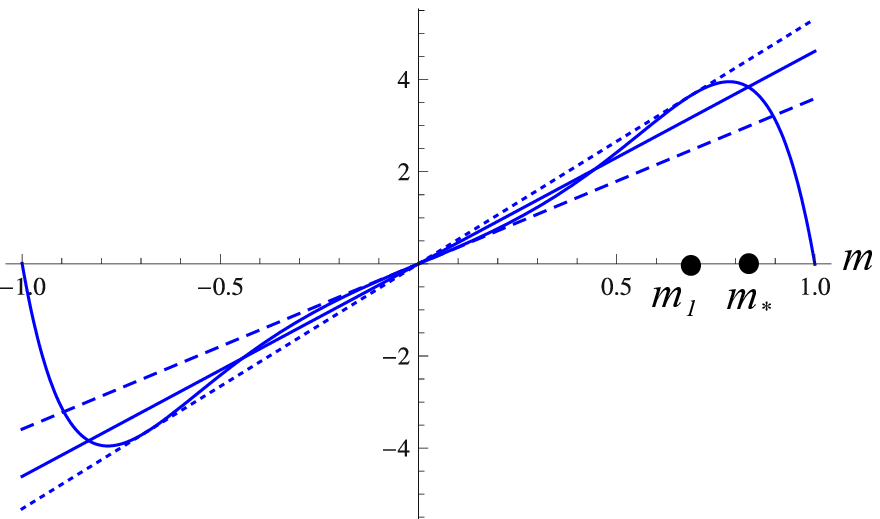}}& 
\subfigure[]{\label{Ct3}\includegraphics[width=0.3 \textwidth]{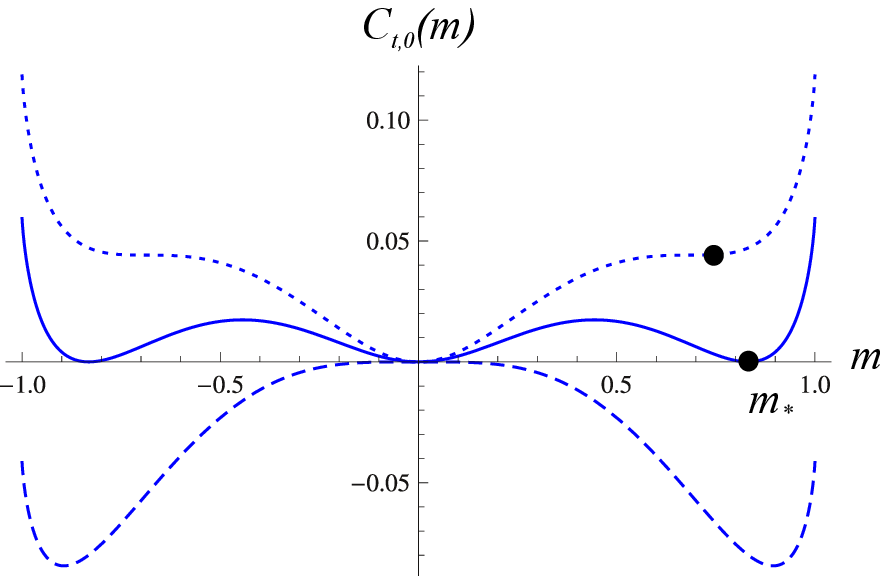}}&
\subfigure[]{\label{CM3}\includegraphics[width=0.3 \textwidth]{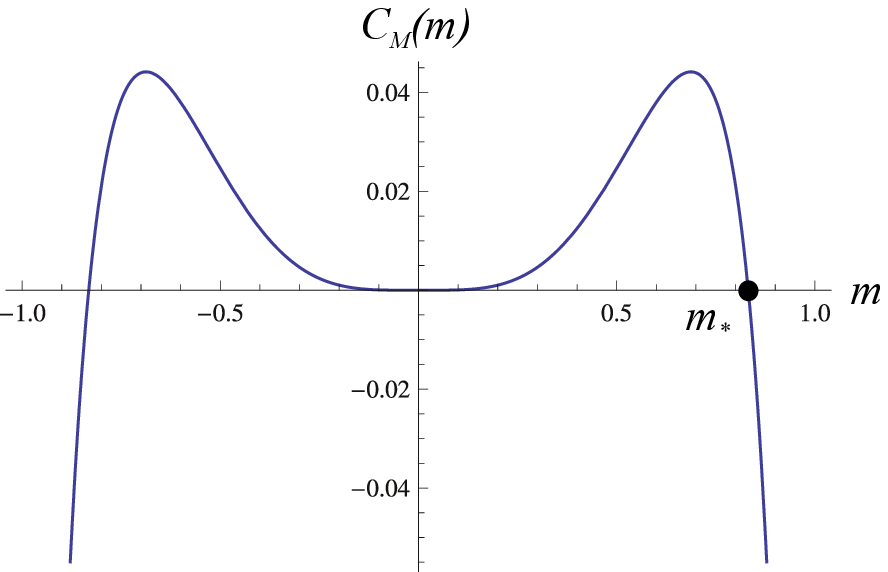}}\\
& Regime (iii), $J=2$. &
\end{tabular}
\end{center}
\caption{(a+d) $m \mapsto a^J(m)$, $m \mapsto l_{t,0}(m)$; (b+e) $m \mapsto 
C_{t,0}(m)$, for $0\leq t<\Psi_c$ (dotted), $t=\Psi_c$ (drawn), $t>\Psi_c$ 
(dashed); (c+f) $m \mapsto C_M(m)$.}
\end{figure}


\subsection{Proof of Theorem~\ref{Cone-Monotonicity}}
\label{S2.4}

\begin{proof}
From \eqref{trajectory} it follows that
\be 
\tilde{m}>m,\,\tilde{t}> t \qquad \Longrightarrow \qquad 
\hat{\phi}^{\tilde{m}}_{\tilde{t},0}(s) > \hat{\phi}^m_{t,0}(s) 
\quad \forall\,0\leq s \leq t.
\ee

\medskip\noindent
(i) First note that, because $C_M(m_*)=0$,
\be
\pd{m_*}{J} = -  \,\frac{\displaystyle\pd{C_M}{J}(m_*)}{\displaystyle\pd{C_M}{m}(m_*)}\;.
\ee
As in Section~\ref{S2.3}, case (iii), we have $(\partial C_M/\partial m)(m_*)<0$. But 
\be
\pd{C_M}{J}(m_*)>0 \qquad \Longleftrightarrow \qquad m_* < \tanh(J m_*),
\ee
which yields the claim because $m_*<m^\infty$.

\medskip\noindent
(ii) The claim is straightforward for $1<J\leq \tfrac32$. For $J>\tfrac32$ we
need to prove that the function $J \to a^J(m_*)/m_*$ is strictly increasing.  In fact,
\be
\begin{split}
\pd{}{J} \left[\frac{a^J(m_*)}{m_*}\right] 
&= \frac{1}{m_*^2}\,\left[- \pd{m_*}{J} \left\{a^J(m_*) - m_* \pd{a^J}{m}(m_*)\right\} 
+ m_*\,\pd{a^J}{J}(m_*) \right]\\[8pt]
&= 1 + \cosh(2Jm_*) - m_*\sinh(2Jm_*) \;.
\end{split}
\ee
The strict positivity of the last expression is equivalent to the inequality 
$m_* < \coth(Jm_*)$, which is satisfied for $J>1$.

\medskip\noindent
(iii) We have $a^J(m) \downarrow a^{\tfrac32}(m)$ as $J \downarrow \tfrac32$ for all 
$m \in (0,1)$, with $a^J$ and $a^{\tfrac32}$ continuous. By Dini's theorem, the 
convergence is uniform.

\medskip\noindent
(iv) Since $\Psi_c(\tilde{J})<\Psi_c(\tfrac32)$, the same argument as in case (iii) 
can be used.

\medskip\noindent
(v) and (vi) are consequences of parts (ii) and (iii) of Theorem \ref{hmfzero}.  
\end{proof}


\section{Proof of Theorem \ref{generaltheorem}}
\label{S3}

In Sections~\ref{S3.1} and \ref{S3.2} we prove that overshoots, respectively,
bifurcations, take place in regime (1), respectively (2)--(3), of Theorem \ref{generaltheorem}. The analysis of the former regime does not distinguish on whether the initial field $h$ is zero or not.


\subsection{Regime (1): Overshoots}
\label{S3.1}

The trick is again to write $t$ as a function of $m$. 
From \eqref{Pequation}, we have
\be
-k^{J,h}(m) + 2m \coth(2t) = m \coth(2t) + \alpha \csch(2t).
\ee
Hence, from \eqref{kbh},
\be
\label{quadratic}
\begin{split}
k^{J,h}(m)& [-k^{J,h}(m) + 2m \coth(2t)]\\
&= [m \coth(2t) - \alpha \csch(2t)]\,[m \coth(2t) + \alpha \csch(2t)]
\end{split}
\ee
which implies
\be	
-k^{J,h}(m)^2 -\alpha^2 = (m^2 - \alpha^2) \coth^2(2t) - 2m k^{J,h}(m) \coth(2t)\;.
\ee
Solving for $t$ we find
\be
\label{tIm}
t_{\rm F}(m) := \left\{\begin{array}{ll}\displaystyle
\frac12 
\acoth\Bigl(\frac{m k^{J,h}(m) + \card\alpha \sqrt{\Phi(m)}}{m^2-\alpha^2}\Bigr)
& \mbox{if } m\neq \pm\alpha\\[10pt]
\displaystyle \frac12 
\acoth\Bigl(\frac{\bigl[k^{J,h}(m)\bigr]^2 + m^2 }{2m \,k^{J,h}(m)}\Bigr)
& \mbox{if } m= \pm\alpha\;,\; m\,k^{J,h}(m)<0\\[10pt]
\displaystyle 0 & \mbox{if } m= \pm\alpha\;,\; m\,k^{J,h}(m)>0
\end{array}\right.
\ee
\be
\label{tDm}
t_{\rm L}(m) := \left\{\begin{array}{ll}\displaystyle
\frac12 
\acoth\Bigl(\frac{m k^{J,h}(m) - \card\alpha \sqrt{\Phi(m)}}{m^2-\alpha^2}\Bigr)
& \mbox{if } m\neq \pm\alpha\\[10pt]
\displaystyle \frac12 
\acoth\Bigl(\frac{\bigl[k^{J,h}(m)\bigr]^2 + m^2 }{2m \,k^{J,h}(m)}\Bigr)
& \mbox{if } m= \pm\alpha\;,\; m\,k^{J,h}(m)>0\\[10pt]
\displaystyle 0 & \mbox{if } m= \pm\alpha\;,\; m\,k^{J,h}(m)<0
\end{array}\right.
\ee
with 
\be
\Phi(m) := k^{J,h}(m)^2 - m^2 + \alpha^2.
\ee
These are times at which $m$ is a stationary point (not necessarily 
a minimum) both for $C_{t,\alpha}(m)$ and for $C_{t,-\alpha}(m)$ [Equation \eqref{quadratic}, and hence the solutions \eqref{tIm} are insensitive to the sign of $\alpha$].  Overshoots and undershoots occur for values of $m$ satisfying (i) $t_F(m)>0$ and $t_L(m)>0$ and (ii) at these times $m$ is a minimum. 

\begin{figure}[htbp]
\begin{center}
\begin{tabular}{cc}
& \multirow{3}{*}{\includegraphics[height=0.5\textwidth]{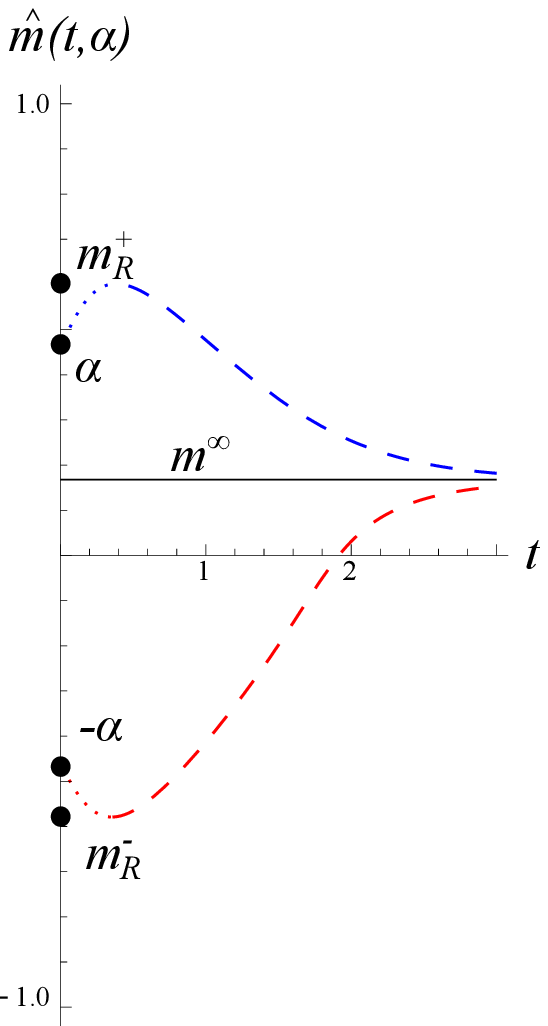}}\\
\includegraphics[width=0.5\textwidth]{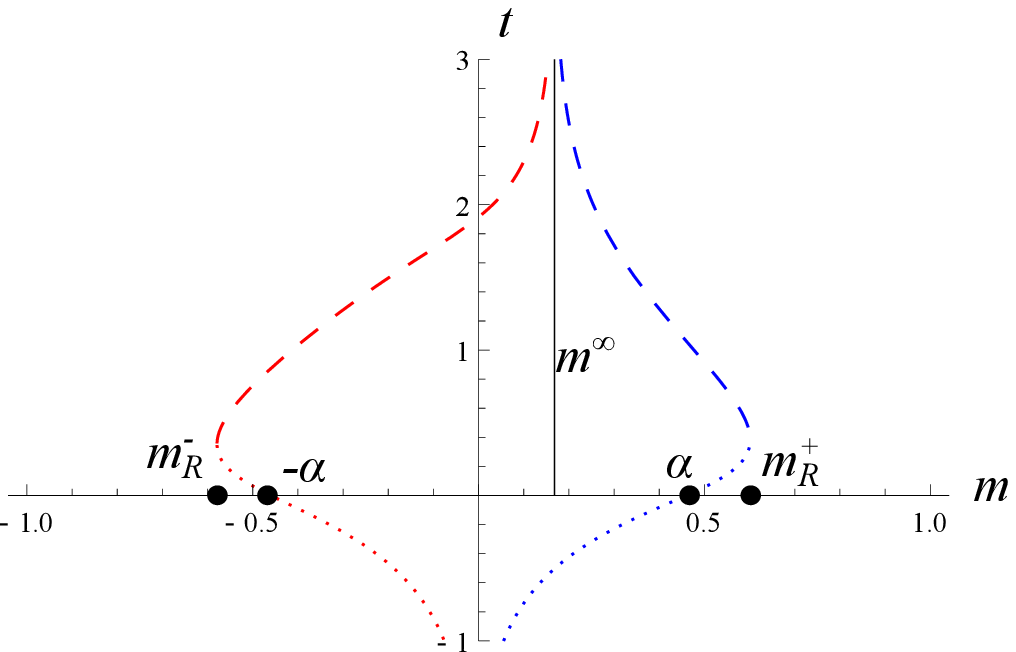} & \\
& \\
& \\
&\\
&\\
&\\
{\small Plot of $t_F$ (dotted), $t_L$ (dashed) for $\alpha$ and $-\alpha$.} &
{\small Plot of $t \mapsto \hat{m}(t)$ for $\alpha$ and $-\alpha$.}
\end{tabular}
\end{center}
\caption{Overshoot for $(J,h,\alpha)$ in Regime (1a) or for $(J,h,-\alpha)$ 
in Regime (1d).
Parameters: $(J,h)=(0.95,0.01), \ \alpha=0.46$.}
\label{Overshoot-R1A-R2B}
\end{figure}
\begin{figure}[htbp]
\begin{center}
\begin{tabular}{cc}
& \multirow{3}{*}{\includegraphics[scale=0.7]{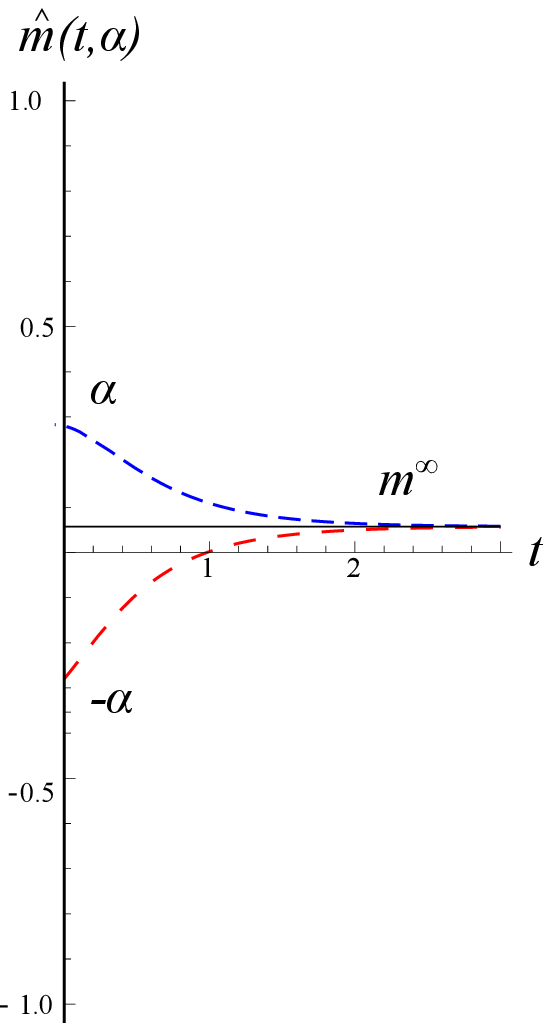}}\\
\includegraphics[scale=0.7]{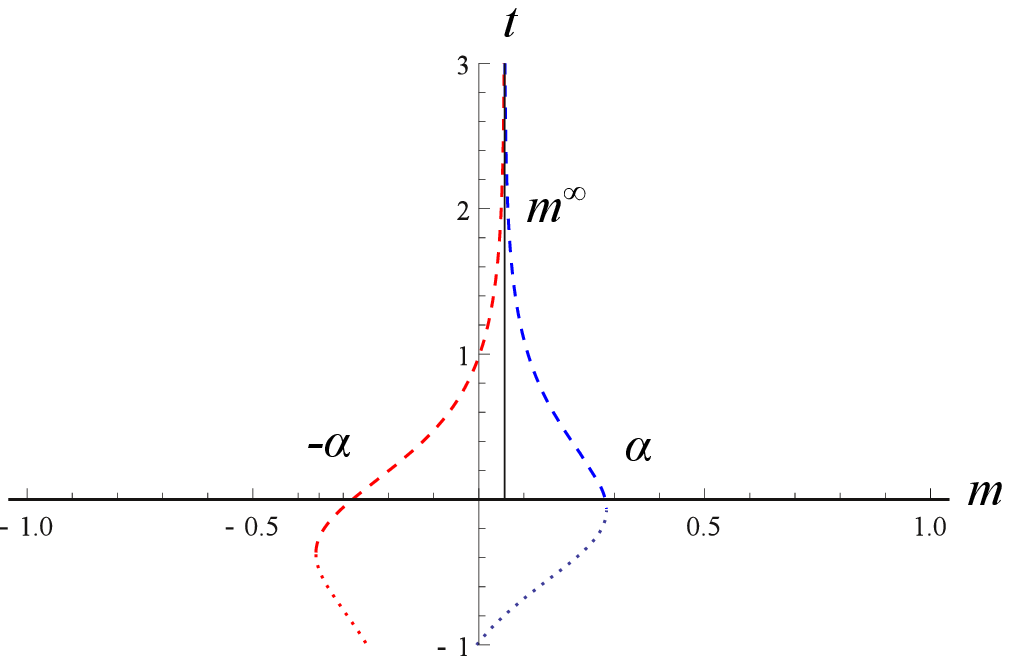} &\\
& \\
&\\
&\\
&\\
{\small Plot of $t_F$ (dotted), $t_L$ (dashed) for $\alpha$ and $-\alpha$.} &
{\small Plot of $t \mapsto \hat{m}(t,\alpha)$ for $\alpha$ and $-\alpha$.}
\end{tabular}
\end{center}
\caption{Absence of overshoot for $(J,h,\alpha)$ in Regime (1b) or 
for $(J,h,-\alpha)$ in Regime (1c). Parameters: $(J,h)=(0.3,0.04), \ \alpha=0.28$.}
\label{NONOvershoot-R1B-R2A}
\end{figure}

\noindent
The $m$-dependence of $t_F$ and $t_L$ is depicted in Figs~\ref{Overshoot-R1A-R2B} and \ref{NONOvershoot-R1B-R2A} for cases in which both functions are injective, i.e., for $\alpha$ 
for which there is only one critical point for each $t$.  In more complicated cases, for instance, 
when overshoot and bifurcation occur simultaneously, there are two or more stationary points only one of which is a minimum.

\noindent
We divide the analysis in four steps. 

\medskip\noindent
{\bf Step A:} \emph{Existence of $z^+$, $z^-$ and $z^{-+}$.}
We observe that there exists a unique $m^\infty>0$ such that $k^{J,h}(m^\infty)=m^\infty$. 
Furthermore,
\be
k^{J,h}(m)=0 \qquad \Longleftrightarrow \qquad \tanh(2h)=\dfrac{-a^J(m)}{b^J(m)} =:A(m).
\ee
This function $A$ has the features depicted in  Figure~\ref{PlotA}: it is odd, satisfies $A(0)=0$, $A(1)=1$ and it is convex with only one global minimum between $0$ 
and $1$. We conclude that there exists a unique $z^+=z^+(J,h)>0$ such that $k^{J,h}(z^+)=0$ and, in addition, 
if $\tanh(2h)<\max\limits_{m \in [-1,0]} A(m)$, then there exist $-1<z^-<z^{-+}<0$ such that
$k^{J,h}(z^-)=k^{J,h}(z^{-+})=0$ (see Fig.~\ref{PlotA}).
	
\begin{figure}[!h]
\begin{center}
\includegraphics[scale=
0.6]{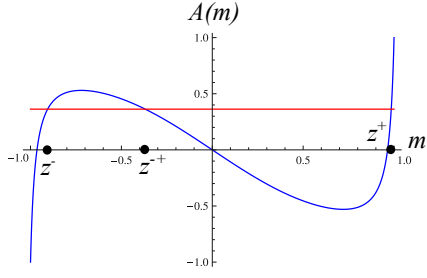}
\end{center}
\caption{Plot of $A$ and intersection with the constant $\tanh(2h)$.}
\label{PlotA}
\end{figure}

\medskip\noindent
{\bf Step B:} \emph{Existence of $m_R^+, \ m_R^-$ and relation between $m_\infty$ and $\alpha$.}

\medskip\noindent
\emph{(Ba) Existence of $m_R^+$:}
$k^{J,h}(\alpha)>0$ together with $\alpha>0$ and Step A imply $\alpha<z^+$. Since 
$\Phi(\alpha)>0$ and $\Phi(z^+)<0$, it follows that there exists a $m_R^+$ such 
that $0<\alpha < m_R^+<z^+$ and 
\begin{equation}
\label{mRplus}
\Phi(m_R^+)=0. 
\end{equation}
The latter in turn implies that $k^{J,h}(m_R^+)^2=m_R^{+^2}-\alpha^2$. This, together with 
$k^{J,h}(m_R^+)>0$, implies that $k^{J,h}(m_R^+)<m_R^+$, which leads to $m^\infty<m_R^+$.
\smallskip\par\noindent
\emph{(Bd) Existence of $m_R^-$}:
As in (Ba), $\Phi(z^-)<0$ and $\Phi(\alpha)>0$ imply that there exists a $m_R^-$ such 
that $z^-<m_R^-<\alpha<0$ and 
\begin{equation}
\label{mRminus}
\Phi(m_R^-)=0.
\end{equation}
\smallskip\par\noindent
\emph{(Bc):} $k^{J,h}(\alpha)>0$ and $\alpha<0$ imply $\alpha< m_\infty$.  This follows from the fact that $k^{J,h}(\alpha)>\alpha$ implies $\alpha< m_\infty$ by \eqref{eq:rr.infty}.
\smallskip\par\noindent
\emph{(Bb):} $k^{J,h}(\alpha)<0$ and $\alpha>0$ imply $\alpha>m_\infty$.  Again, this is a consequence of  \eqref{eq:rr.infty}.

\medskip\noindent
{\bf Step C:} \emph{Consequence of the positivity of times.}
Only positive solutions of Equation \eqref{tIm} are of interest. This implies the constraints
\be
\label{eF}
t_F(m)> 0 \quad \Longleftrightarrow \quad 
\eta_F(m) := m\, k^{J,h}(
m)+ (\alpha^2-m^2)+ |\alpha| \sqrt{\Phi(m)}
\begin{cases}
>0 &\text{if } m^2> \alpha^2,\\
<0 &\text{if } m^2<\alpha^2,
\end{cases}
\ee
and
\be
\label{eL}
t_L(m)> 0 \quad \Longleftrightarrow \quad 
\eta_L(m):=m\, k^{J,h}(m) + (\alpha^2-m^2)- |\alpha| \sqrt{\Phi(m)}
\begin{cases}
>0 &\text{if } m^2> \alpha^2,\\
<0 &\text{if } m^2<\alpha^2.
\end{cases}
\ee
The functions $\eta_F$ and $\eta_L$ satisfy
\be
\eta_F(\alpha) = \alpha k^{J,h}(\alpha) + |\alpha|  |k^{J,h}(\alpha)|=
\begin{cases}
>0 &\text{if } \alpha \ k^{J,h}(\alpha)>0,\\
=0 &\text{if } \alpha \ k^{J,h}(\alpha)<0,
\end{cases}
\ee
\be
\eta_L(\alpha) = \alpha k^{J,h}(\alpha) - |\alpha| |k^{J,h}(\alpha)|=
\begin{cases}
=0 &\text{if } \alpha \ k^{J,h}(\alpha)>0,\\
<0 &\text{if } \alpha \ k^{J,h}(\alpha)<0.
\end{cases}
\ee
Also, from \eqref{eq:rr.infty},
\begin{equation}
\label{propertiese}
\begin{array}{lcl}
\eta_F(m^\infty) = 2|\alpha|^2, &\qquad& \eta_F'(m^\infty)= 2 m^\infty {k^{J,h}}'(m^\infty),\\
\eta_L(m^\infty) = 0, &\qquad& \eta_L'(m^\infty)=0.
\end{array}
\end{equation}
Last line implies that $m^\infty$ is a root of $\eta_L$ but there is no change of sign around it. Finally, from expressions \eqref{eF}-\eqref{propertiese} we conclude that:
\begin{itemize}
\item The zeros of $\eta_F, \ \eta_L$ are a subset of $\{m^\infty, \pm \alpha\}$.
\item  The intervals in which $\eta_F$ and $\eta_L$ satisfies the constrains \eqref{eF}-\eqref{eL} are:
\end{itemize}
\begin{center}
\begin{tabular}{|l|c|c|c|c|}
\hline
Regime & $(1a)$ & $(1b)$ & $(1c)$ & $(1d)$ \\
\hline
Condition& $\alpha>0, \ k^{J,h}(\alpha)>0$ & $\alpha>0, \ k^{J,h}(\alpha)<0$& $\alpha<0, \ k^{J,h}(\alpha)>0$ & $\alpha<0, \ k^{J,h}(\alpha)<0$\\
\hline
$t_F>0$ & $[\alpha,m_R^+]$ & $\emptyset$ & $\emptyset$ & $[m_R^-,\alpha]$ \\
\hline
$t_L>0$ & $[m^\infty,m_R^+]$ & $[m^\infty,\alpha]$ & $[\alpha,m^\infty]$ & $[m_R^-,m^\infty]$ \\ \hline
\end{tabular}
\end{center}

We observe that in regime (1a) each value of $m\in[\alpha\wedge m_\infty, m_R^+]$ is attained at two different times $t_F$ (``First") and $t_L$ (``Last").  The same happens in regime (1d) for
$m\in[m_R-,\alpha\vee m_\infty]$.  These phenomena correspond respectively, to an over and an under shoot.  The proof is completed by showing that the trajectories have the right monotonicity properties.

\begin{figure}[hbtp]
\begin{center}
\includegraphics[width=0.45 \textwidth]{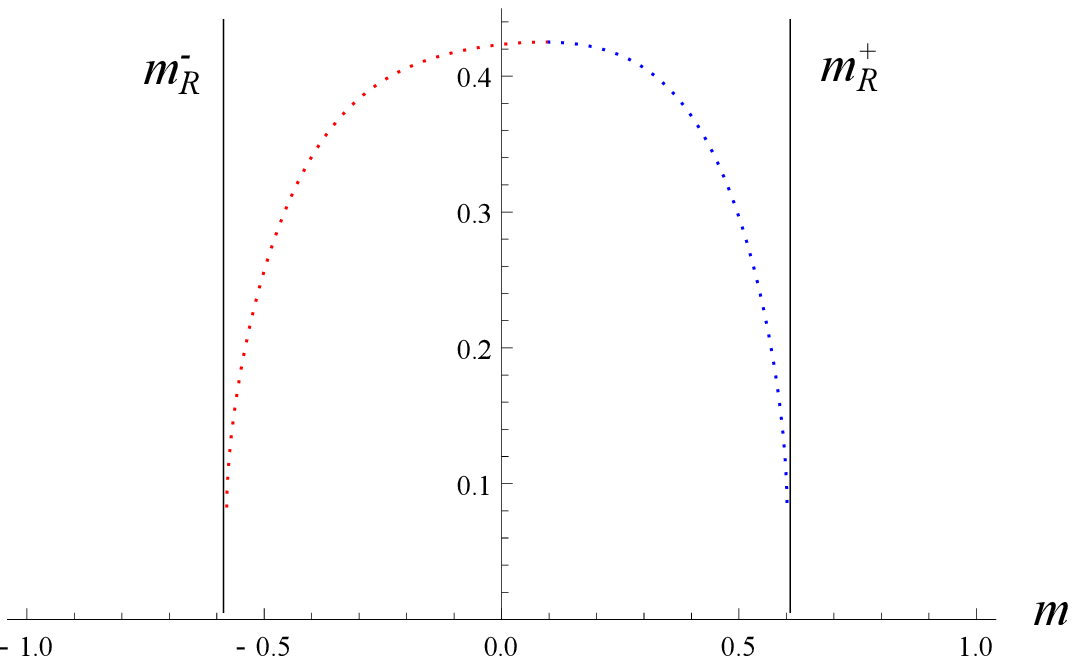}
\includegraphics[width=0.45 \textwidth]{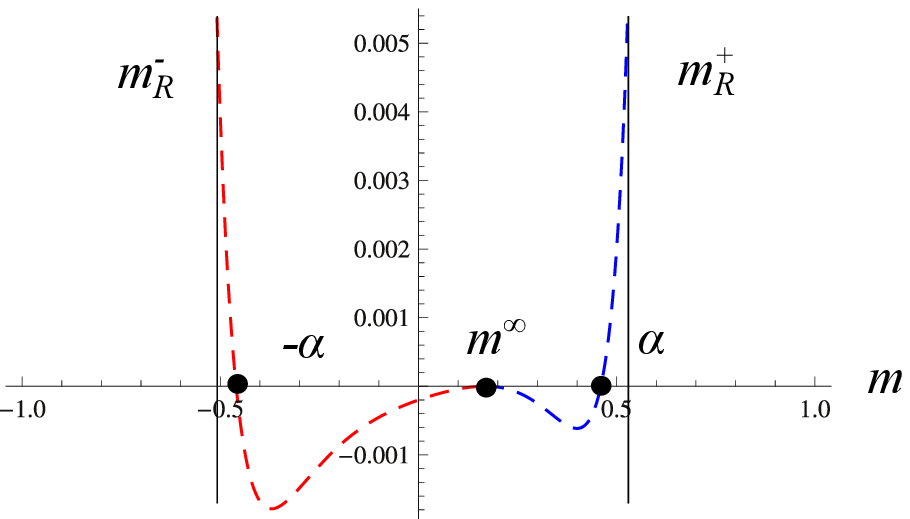}\\
\end{center}
\caption{$\eta_F$ (dotted) and $\eta_L$ (dashed) for $(J,h,\alpha)$ in regime 1(a) and for $(J,h,-\alpha)$ in regime 1(d). Parameters: $(J,h)=(0.95,0.01), \ \alpha=0.46$ as in Fig. \ref{Overshoot-R1A-R2B}.}
\label{eDeI}
\end{figure}
\medskip\noindent
{\bf Step D:} \emph{Monotonicity.}\\
By using implicit derivation we get  
\be\label{eq:mono}
\pd{\hat{m}}{t}(t)=\frac{\pd{\l}{t}(\hat{m})}{[\k]'(\hat{m})-\l'(\hat{m})}=\frac{2 \csch(2t) \Big\{ \hat{m} \csch(2t)-\alpha \coth(2t) \Big\} }{[\k]'(\hat{m})-\coth(2t)}\;.
\ee
On the other hand, if $\hat{m}=\hat{m}(t)$ is a critical point $\Big(\k(\hat{m})=\l(\hat{m})\Big)$, then $\Phi(\hat{m})=(\hat{m} \csch(2t)-\alpha \coth(2t))^2$. Hence,

$$\pd{\hat{m}}{t}(t)=0 \quad \Longleftrightarrow \quad \Phi(\hat{m}(t))=0.$$
Splitting in cases according to different values of $\hat{m}$, it is not hard to conclude from \eqref{mRplus} and \eqref{eq:mono} the following monotonicity properties of the trajectories. 

\begin{center}
\begin{tabular}{|l|c|c|c|c|}
\hline
Regime & $(1a)$ & $(1b)$ & $(1c)$ & $(1d)$ \\
\hline
$\hat{m}$ incr. & $0<t<t(m_R^+)$ & $\emptyset$ & $0<t<\infty$ & $t(m_R^-)<t<\infty$ \\
\hline
$\hat{m}$ decr. & $t(m_R^+)<t<\infty$ & $0<t<\infty$ & $\emptyset$ & $0<t<t(m_R^-)$ \\ \hline
\end{tabular}
\end{center}

\noindent
This concludes the proof of part 1 of Theorem \ref{generaltheorem}.

\subsection{Bifurcation}
\label{S3.2}

Bifurcation proofs rely on the following facts.  
\begin{itemize}
\item[\bf (B1)]
For short times there is a unique critical point, close to $\alpha$.

\item[\bf (B2)] 
Therefore, in order for bifurcation to occur a local maximum 
and a local minimum must appear in the course of time. Given condition \eqref{Pequation}, usual arguments imply that two (or more) stationary points appear at times larger than $\tilde t$ if the curves $l_{\tilde{t},\alpha}$ and $k^{J,h}$ become tangent at a certain magnetization $\tilde m$. The pairs $(\tilde{m},\tilde{t})$ are determined by the following two equations 
(a similar argument was used in the proof of Theorem~\ref{hmfzero}(iii)):	

\begin{equation}
\label{tagencysystem}
\begin{split}
\bigl[{k^{J,h}}\bigr]'(\tilde{m}) &= \coth(2\tilde{t}),\\
k^{J,h}(\tilde{m}) &= \tilde{m} \coth(2\tilde{t})-\alpha \csch(2\tilde{t}).  
\end{split}
\end{equation}
Inserting the first equation into the second, we get
\begin{equation}
\label{tangency}
F(\tilde{m}) := \dfrac{\tilde{m} \bigl[{k^{J,h}}\bigr]'(\tilde{m})-k^{J,h}(\tilde{m})}
{\csch\bigl[\acoth\bigl(\bigl[{k^{J,h}}\bigr]'(\tilde{m})\bigr)\bigr]} = \alpha.
\end{equation}
We are left with the task of determining whether or not this equation has solutions.
Note that 
\be
\label{criticalF}
F'(m)=0 \quad \Longleftrightarrow \quad \bigl[{k^{J,h}}\bigr]''(m)=0 
\quad \text{or} \quad m=k^{J,h}(m)\bigl[{k^{J,h}}\bigr]'(m).
\ee
In what follows all the assertions about $F$ can be checked by using the equivalence 
in \eqref{criticalF} and doing a straightforward analysis of $k^{J,h}$. 

\item[\bf (B3)] $t \mapsto C_{t,\alpha}$ is continuous with respect to $\|\cdot\|_\infty$. 
Hence, when a new minimum appears it cannot be a global one.
Both $\alpha \mapsto C_{t, \alpha}$ and $h \mapsto C_{t,\alpha}$ are also continuous.

\item[\bf (B4)] When $t \to \infty$ we have two global minima $\pm m^\infty$ if $h=0$. 
If $h>0$ the symmetry is broken and there is only one global minimum $m^\infty>0$.
\end{itemize}

Whenever a local 
maximum/local minimum appears (disappears), we will refer to this behaviour as 
LMLMA (LMLMD).
We proceed by looking at $h=0$ and $h\neq 0$ separately.

\medskip\noindent
{\bf Part (2)} ($h=0$, see Fig.~\ref{Fhzero}). The scenario for $\alpha=0$ has already been proven in Theorem \ref{hmfzero}.  We concentrate on $\alpha>0$; this is no loss of generality due to the antisymmetry of 
$F$. 

\emph{Claim:} Whenever $\alpha>0$, negative solutions of \eqref{tangency} can not cause  bifurcations. 
In fact, let $t_{nc}$ be the time at which the critical point in the negative side emerge.
Let also, 
$$d(t):=C_{t,\alpha}(\overline{m}_{-}(t,\alpha))-C_{t,\alpha}(\hat{m}(t,\alpha)), \qquad t \geq t_{nc},$$
where $\overline{m}_{-}(t,\alpha)$ is the negative \emph{local} [because of {\bf (B3)}] minimum  and $\hat{m}(t,\alpha)$ is the global minimum of $C_{t,\alpha}$. The last one is positive due to {\bf (B1)}. and the supposition of $\alpha>0$.
\noindent
By definition, $d(t_{nc})>0$ and by {\bf (B4)} $\lim_{t \to \infty} d(t)=0$.
Doing calculations similar to \eqref{dC/dt} and using that $\overline{m}_{-}(t,\alpha)$, $\hat{m}(t,\alpha)$ are both critical points, we get that $d'(t)<0 \ \forall \ t > t_{nc}$. This proves the claim.

In what follows we focus in equation \eqref{tangency}. Owing to the previous claim, in order to bifurcation to occur a positive solution of \eqref{tangency} is needed.
\begin{itemize}
\item[(2a-b)]  
If $0<J\leq \tfrac{3}{2}$, then $F'(m)<0$ for all $m \in (-1,+1)$. Hence, for all 
$\alpha>0$ there is only one solution of \eqref{tangency}. This solution turns out to be negative and hence it can not correspond to a bifurcation.  
\item[(2c)]  
If $J > \tfrac32$, $F$ has only one global maximum on the positive side, with value $U_B=U_B(J)>0$. 
Combining {\bf B1.}-{\bf B4.}, we get that there is bifurcation if and only if $\alpha \in[0,U_B]=\mathrm{Im}(F|_{[0,1]})$. 
\end{itemize}

\medskip

\begin{figure}[h]
\begin{center}
\begin{tabular}{cc}
\includegraphics[width=0.45\textwidth]{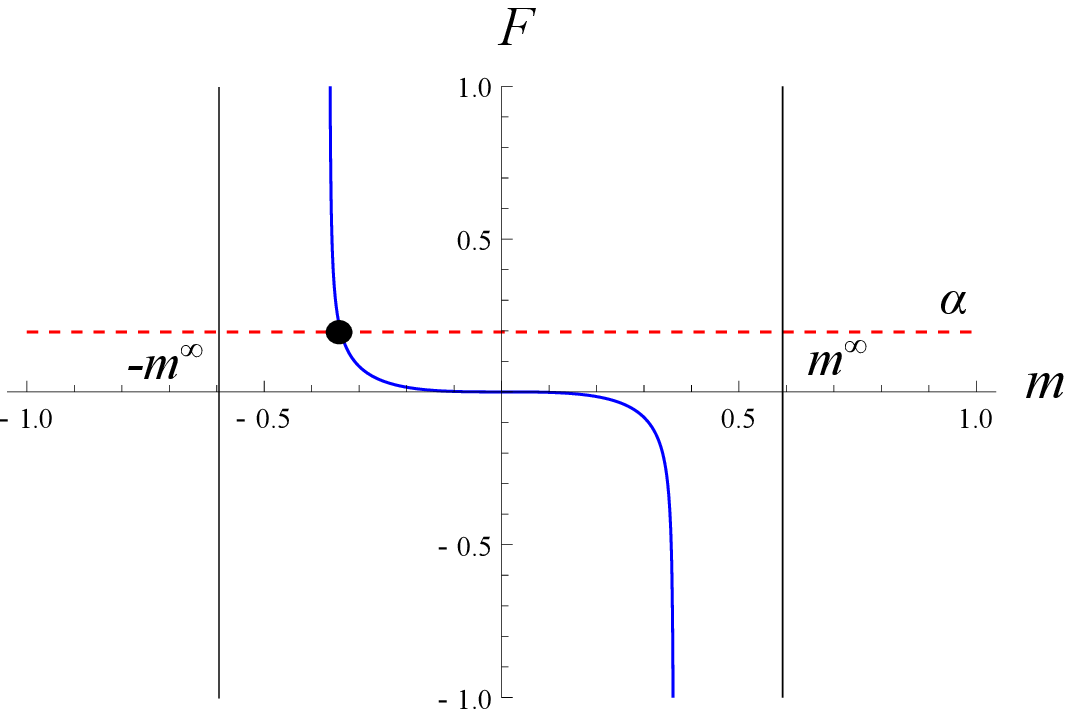} 
&\includegraphics[width=0.45\textwidth]{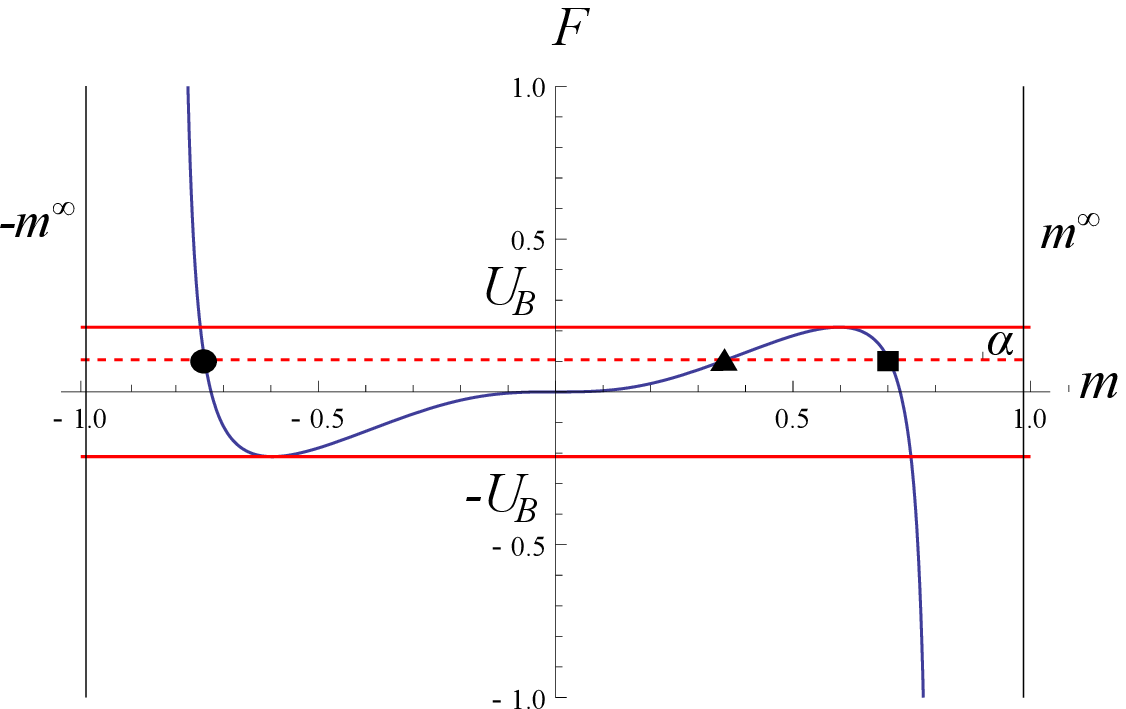} \\
\includegraphics[width=0.45\textwidth]{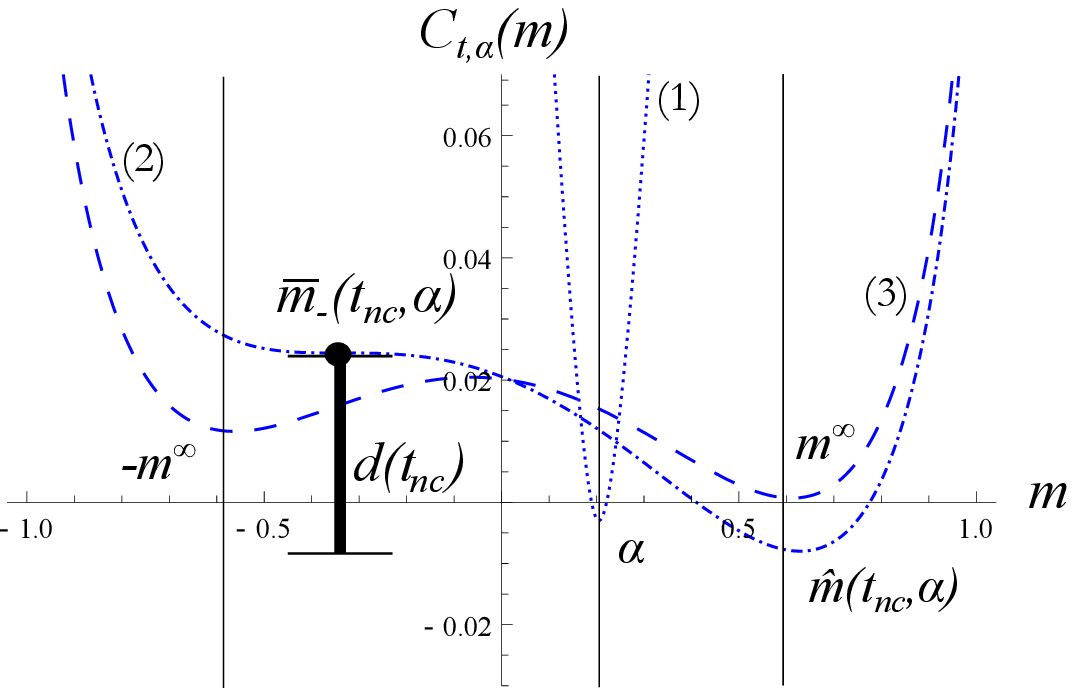} 
& \includegraphics[width=0.45\textwidth]{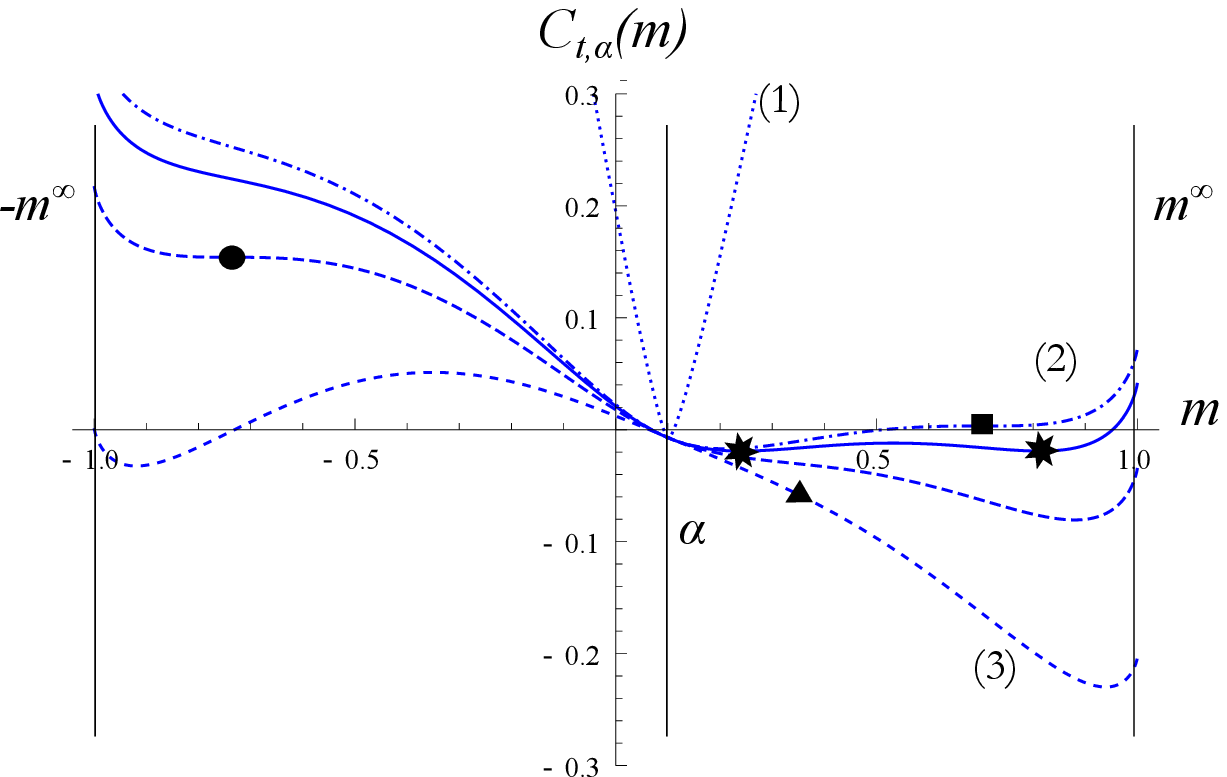} \\
Regime (2a-b), $J=1.15$, $h=0$ & Regime (2c), $J=2.5$, $h=0$. \\
\end{tabular}
\end{center}
\caption{First row: Plot of $m \mapsto F(m)$. $\newmoon=$ LMLMA not leading to 
bifurcation; $\blacktriangle=$ LMLMD; $ \blacksquare=$ LMLMA leading to 
bifurcation; $\bigstar=$ bifurcation . Second row: Plot of $m \mapsto C_{t,\alpha}(m)$ for different times. 
Short time (dotted) bifurcation (solid), long time (dashed).}
\label{Fhzero}
\end{figure}

\medskip\noindent
{\bf Part (3)} ($h>0$, see Fig.~\ref{Fhnonzero}).\\
\emph{Remark:} 
If $h>0$, then {\bf B4.} (Symmetry breaking) allows the appearance of solutions of \eqref{tangency} leading to bifurcations.

Once more, let study the different scenarios for $F$ when $h>0$.
\begin{itemize}
\item[(3a)] 
If $0<J \leq 1$, then $\mathrm{Im}(F)\subset [-1,+1]^c$. Therefore \eqref{tangency} has no solution for any $\card\alpha\le 1$ and there is no bifurcation.
\item[(3b)]
If $1<J \leq \tfrac32$, then $F$ has a unique maximum for $m\in[0,1]$ with value 
$U_B=U_B(J,h)<0$, and $[-1,U_B]=\mathrm{Im}(F|_{[0,1]})$. 
Arguing as in the claim of Part 2 for $\alpha>U_B$, we conclude that there is bifurcation if and only if $\alpha \in [-1,U_B]$. 

\item[(3c)] Assume $J> \tfrac32$.
\begin{enumerate}
\item 
For $h>0$ small enough the behaviour is ``close'' to the $h=0$ case due to the continuity of $h \mapsto C_{t,\alpha}$ with respect to the infinite norm.\\ \noindent
Indeed, there exists $L_B:=\min_{[-1,0]}F \approx -U_B(J,0)$ and $U_B:=\max_{[0,1]}F \approx U_B(J,0)$, with $(-1,U_B]=\mathrm{Im}(F|_{[0,1]})$. There are different regimes for $\alpha$:
\begin{enumerate}
\item 
For $\alpha<L_B$, there is a unique solution of \eqref{tangency}, which is in the positive side, leading to a bifurcation [because of {\bf (B4)} ]. 
\item 
For $0>\alpha \gtrapprox L_B$ there are both a negative and a positive solution to \eqref{tangency}.  
Both lead to bifurcations, the negative one by continuity {\bf B3.} and the positive due to {\bf B3.}-{\bf B4.}.
The negative solution appears earlier in time.  We write $s_B=s_B(\alpha)$ for the time of the first (negative side)
bifurcation and $t_B=t_B(\alpha)$ for the second bifurcation time.
\item 
For $0<\alpha \lessapprox U_B$ there is LMLMA on the positive and negative sides.
As in the $h=0$ case, the negative one does not lead to a bifurcation, and thus only one bifurcation occurs, which happens to be in the positive side.

\item 
By the continuity property {\bf (B3)} and the monotonicity of $\alpha \mapsto s_B(\alpha)$ (proved below) , the two previous regimes coalesce, leading to an intermediate value $M_T\in(L_B,U_B)$ such that trifurcation occurs at $\alpha=M_T$.
\item 
For $\alpha >U_B$ there is no positive solution to \eqref{tangency} with $\alpha\in[0,1]$.  Hence, no bifurcation occurs.
\end{enumerate}
\item The limit $h \to \infty$ in \eqref{tangency} yields
$$\frac{m \Big({a^J}'(m)+{b^J}'(m)\Big)-\Big(a^J(m)+b^J(m)\Big)}{\Big({a^J}'(m)+{b^J}'(m)\Big)}=\alpha.$$
Hence, for $h>0$ large enough we get a behaviour similar to (3b), but with $U_B(J,h)>0$. 

\item The existence of $h^*$ follows from the continuity of the function $h \mapsto C_{t,\alpha}$  commented in {\bf (B3)} with respect to $h$.
\end{enumerate}

\end{itemize}

\begin{figure}[hbtp]
\begin{center}
\begin{tabular}{cc}
\includegraphics[width=0.4\textwidth]{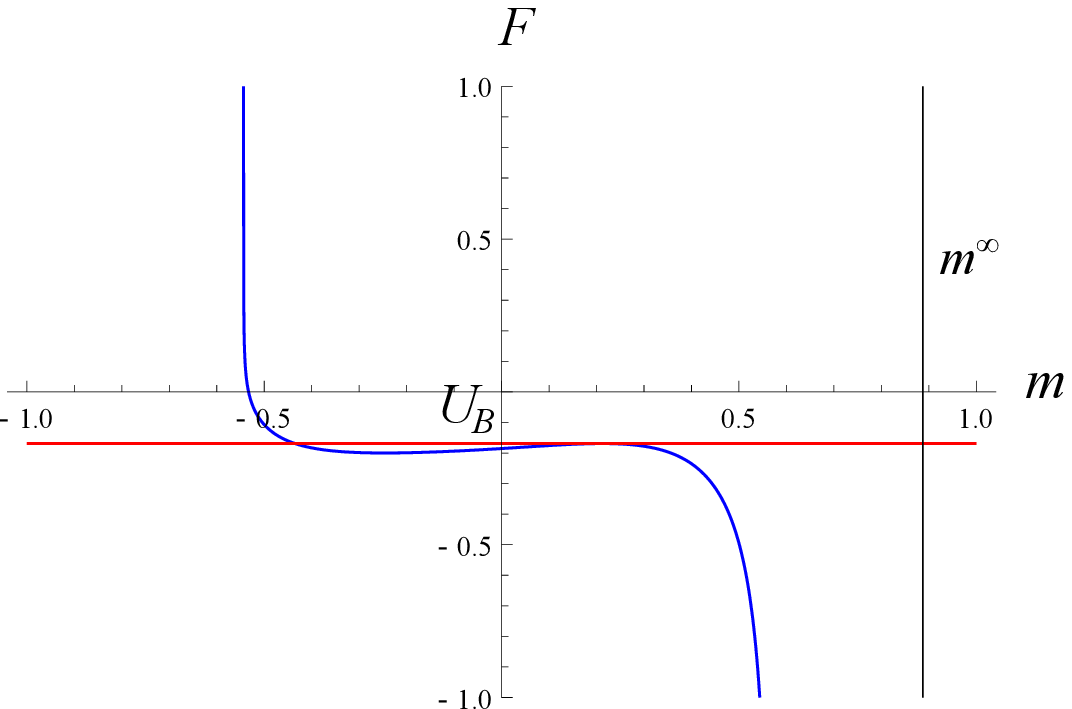} & 
\includegraphics[width=0.4\textwidth]{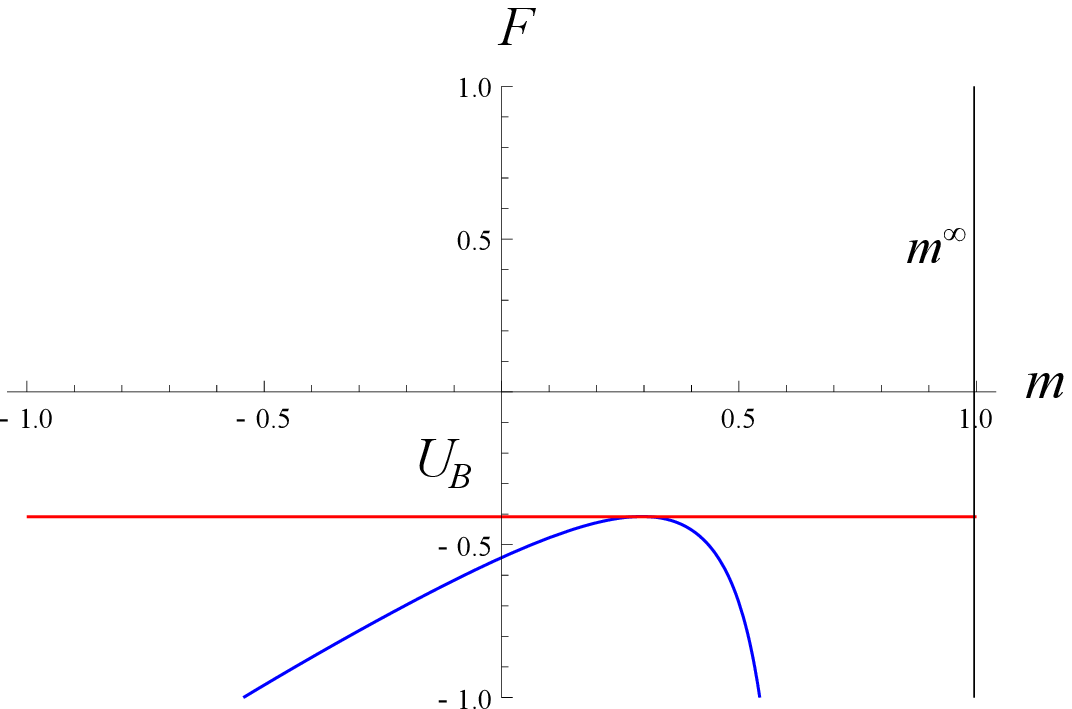}  \\
(3b), $h$ small, \quad $J=1.42$, $h=0.15$ & (3b), $h$ large,\quad $J=1.42$, $h=1.6$\\
\includegraphics[width=0.4\textwidth]{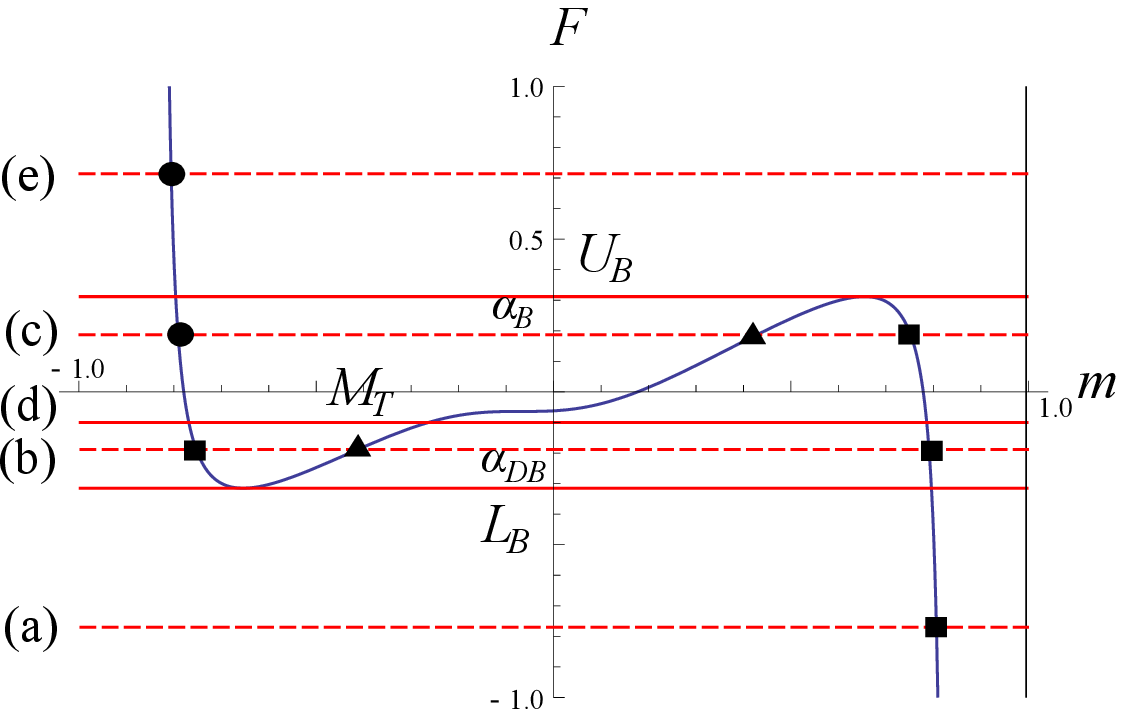} 
& \includegraphics[width=0.4\textwidth]{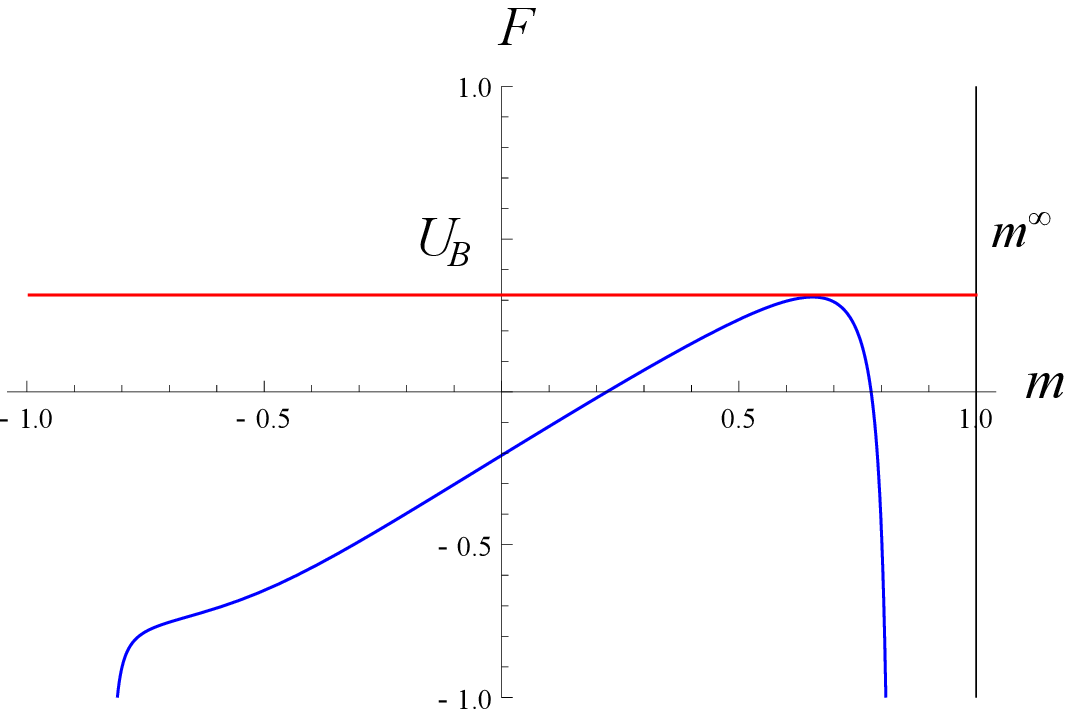} \\
(3c), $h$ small, \quad $J=2.9$, $h=0.15$ &(3c), $h$ large, \quad $J=2.9$, $h=1.6$
\end{tabular}
\end{center}
\caption{Plot of $m \mapsto F(m)$ for different regimes of $J$ when $h>0$.
$\newmoon=$ LMLMA not causing a bifurcation; $\blacktriangle=$ LMLMD; 
$\blacksquare=$ LMLMA causing a bifurcation.}
\label{Fhnonzero}
\end{figure}


\subsubsection{Monotonicity of the functions $t_B(\alpha)$ and $s_B(\alpha)$}
\label{S2.3.3}

The bifurcation times are characterized by the following equations:
\be
\begin{array}{lll}
k^{J,h}(\hat{m}_1) &=& l_{t_B,\alpha}(\hat{m}_1),\\
k^{J,h}(\hat{m}_2) &=& l_{t_B,\alpha}(\hat{m}_2),\\
C_{t_B,\alpha}(\hat{m}_1) &=& C_{t_B,\alpha}(\hat{m}_2).
\end{array}
\ee
The first two equations say that $\hat{m}_1$ and $\hat{m}_2$ are stationary points 
at the same time $t_B$, while the third one establishes the equality of costs at this time $t_B$.  Taking the derivative with respect to $\alpha$ 
of the third equation we get
\be
\pd{t_B}{\alpha} = -\dfrac{\displaystyle \pd{C_{t_B,\alpha}}{\alpha}(\hat{m}_2)-\pd{C_{t_B,\alpha}}{\alpha}(\hat{m}_1)}
{\displaystyle\pd{C_{t_B,\alpha}}{t_B}(\hat{m}_2)-\pd{C_{t_B,\alpha}}{t_B}(\hat{m}_1)}\;.
\ee
A straightforward computation using the first two equations shows that $\pd{t_B}{\alpha}<0$, which implies
that $\alpha\mapsto t_B(\alpha)$ is continuous and decreasing. A similar argument shows that $\alpha\mapsto
s_B(\alpha)$ is continuous and increasing.


\end{document}